\documentclass{article}
\usepackage{amsfonts}\usepackage{amssymb}
\usepackage{graphicx}\usepackage{amsmath}
\usepackage{abstract}
\usepackage{color}
\usepackage{graphics}   \usepackage{amsthm}

\numberwithin{equation}{section}
\providecommand{\keywords}[1]{\textbf{Key words} #1}

\title{Non-conforming  Crouzeix-Raviart element approximation for Stekloff eigenvalues in inverse scattering }

\author{Yidu Yang, Yu Zhang, Hai Bi
 \\\\
{\small School of Mathematical Sciences,
Guizhou Normal University,}\\{\small  Guiyang,  $550001$,  China}
}

\begin{document}
\date{}
\maketitle

\begin{abstract}
In this paper, we use the non-conforming Crouzeix-Raviart element method to solve a Stekloff eigenvalue problem arising in inverse scattering. The weak formulation corresponding to this problem is
non-selfadjoint and does not satisfy $H^{1}$-elliptic condition,
and its Crouzeix-Raviart element discretization does not meet the Strang lemma condition.
We use the standard duality techniques to prove an extension of Strang lemma.
And we prove the convergence and error estimate of discrete eigenvalues and eigenfunctions using the spectral perturbation theory for compact operators.
Finally, we present some numerical examples
not only on uniform meshes but also in an adaptive refined meshes to show that the Crouzeix-Raviart method is efficient for computing real and
complex eigenvalues as expected.
\end{abstract}


\section{Introduction}

Steklov eigenvalue problems have important physical
background and many applications. For instance, they appear
in the analysis of stability of mechanical oscillators immersed
in a viscous fluid (see \cite{conca} and the references therein), in the study
of surface waves \cite{bergman}, in the study of the vibration modes
of a structure in contact with an incompressible fluid
\cite{bermudez} and in the analysis of the antiplane shearing on a system
of collinear
faults under slip-dependent friction law \cite{bucur}. Hence, the finite element methods for solving these problems have attracted more and more scholars' attention.
Till now, systematical and profound studies on the finite element approximation mainly
focus on Steklov eigenvalue problems
which satisfy $H^{1}$-elliptic condition (see, e.g.,
\cite{alonso,andreev,armentano1,armentano2,bermudez,bi,bramble,cao,garau,liq2,lim,russo,xie,yang1} and the references therein).\\
Recently Cakoni et al. \cite{cakoni1} study a new Stekloff eigenvalue problem arising from the inverse scattering theory:
 \begin{eqnarray}\label{s1.1}
\Delta u +k^{2} n(x)u=0~~in~ \Omega,~~~~
\frac{\partial u}{\partial \gamma}=-\lambda u~~ on~ \partial\Omega,
\end{eqnarray}
where $\Omega$ is a bounded domain in $\mathbb{R}^{d}$ ($d=2,3$),
$\frac{\partial u}{\partial \gamma}$ is the outward normal derivative,
$k$ is the wavenumber and
$
n(x) = n_{1}(x)+i\frac{n_{2}(x)}{k}
$
is the index of refraction that
is a bounded complex valued function with
$n_{1}(x)>0$ and $n_{2}(x)\geq 0$.\\
\indent Note that the weak formulation of (\ref{s1.1}) (see (\ref{s2.1})) does not satisfy $H^{1}$-elliptic condition.
Cakoni et al. \cite{cakoni1}
 analyze the mathematics properties of (\ref{s1.1}) and use conforming finite element method to solve it.
Liu et al. \cite{liu} then study error estimates of conforming finite element eigenvalues for (\ref{s1.1}).
\\
\indent In this paper, we will study the non-conforming Crouzeix-Raviart element
(C-R element) approximation for the problem (\ref{s1.1}).
The C-R element was first introduced by  Crouzeix and Raviart in \cite{crouzeix} in 1973
to solve the stationary Stokes equation. It was also used to solve linear elasticity equations (see \cite{falk,brenner2}),
the Laplace equation/eigenvalues (see \cite{armentano3,boffi,brenner,carstensen5,carstensen6,carstensen1,duran}), Darcy's equation \cite{ainsworth},
Steklov eigenvalue (see\cite{alonso,bi,liq2,russo,yang1}) and so on.
The features of our work are as follows:
\begin{enumerate}
\item
As we know, the convergence and error estimates of the  non-conforming finite element method for an eigenvalue problem is
based on the convergence and error estimates of the  non-conforming finite element method for the corresponding source problem, and
Strang lemma (see \cite{strang}) is a fundamental analysis tool.
However, the sesquilinear form in the C-R element discretization here does not meet the
Strang lemma condition.
To overcome this difficulty, referring to \S 5.7 in \cite{brenner}, we use the standard duality techniques to prove an extension of Strang lemma (see Theorem 2).
Based on the theorem, we prove the convergence and error estimates of the C-R method for the corresponding source problem.
The current paper, to our knowledge, is the first investigation of applying and extending Strang lemma to elliptic boundary value problem that the
corresponding sesquilinear form is non-selfadjoint and not coercive.
\item
Cakoni et al. \cite{cakoni1} write (\ref{s1.1}) as an equivalent eigenvalue problem of the Neumann-to-Dirichlet operator $T$. In this paper,
we write the C-R element approximation of (\ref{s1.1}) as an equivalent eigenvalue problem of the discrete operator $T_{h}$, and prove $T_{h}$ converges $T$ in the sense of norm in $L^2(\partial\Omega)$, thus using Babuska-Osborn spectral approximation
theory \cite{babuska} we prove first the convergence and error estimates of C-R finite
element eigenvalues and eigenfunctions for the problem (\ref{s1.1}).
\item
We implement some numerical experiments
not only on uniform meshes but also in adaptive refined meshes. It can be seen that the C-R
method is efficient for computing real and complex eigenvalues as expected.
In addition, we discover, when the index of refraction $n(x)$ is real and
$\Omega$ is the
L-shaped domain or
the square with a slit, the C-R element eigenvalues approximate the exact ones from above,
 and numerical results in \cite{cakoni1,liu} show conforming finite element eigenvalues approximate the exact ones from below, thus we get the upper and lower bounds of eigenvalues.
\end{enumerate}

\indent It should be pointed out that the theoretical analysis and conclusions in this paper are also valid for the extension Crouzeix-Raviart element \cite{huj}.\\
 \indent In
this paper, regarding the basic theory of  finite
element methods, we refer to \cite{babuska,brenner,chatelin,ciarlet,oden,shiwang}.\\
\indent Throughout this paper, $C$ denotes a positive constant
independent of $h$, which may not be the same constant in different
places. For simplicity, we use the symbol $a\lesssim b$ to mean that
$a\leq C b$.

\section{Preliminary}

\indent In this paper, we assume $\Omega \subset \mathbb{R}^{d}$ ($d=2,3$) is a polygonal ($d=2$) or polyhedral ($d=3$) domain,
and assume $n\in L^{\infty}(\Omega)$.
 Let $H^{\rho}(\Omega)$ denote the Sobolev
space on $\Omega$
with the norm $\|\cdot\|_{\rho,\Omega}$ (denoted by $\|\cdot\|_{\rho}$ for simplicity) and the seminorm $|\cdot|_{\rho,\Omega}$ (denoted by $|\cdot|_{\rho}$ for simplicity) and $H^{0}(\Omega)=L^{2}(\Omega)$, and let
$H^{\rho}(\partial\Omega)$ denotes the Sobolev space on $\partial\Omega$ with the norm $\|\cdot\|_{\rho,\partial\Omega}$ and the seminorm $|\cdot|_{\rho,\partial\Omega}$.\\
\indent  Cakoni et al. \cite{cakoni1} give the weak form of (\ref{s1.1}): Find $\lambda\in
\mathbb{C}$, $u \in H^{1}(\Omega)\setminus\{0\}$, such that
\begin{eqnarray}\label{s2.1}
a(u,v)=-\lambda <u,v>,~~~\forall v\in H^{1}(\Omega),
\end{eqnarray}
where
\begin{eqnarray*}
a(u,v)=(\nabla u, \nabla v)-(k^2nu, v),~~~(u,v)=\int \limits_{\Omega}u\bar{v}dx,~~~<u,v>=\int \limits_{\partial\Omega}u\bar{v}ds.
\end{eqnarray*}
The source problem associated with (\ref{s2.1}) is as follows:
Find $\varphi\in H^{1}(\Omega)$ such that
\begin{eqnarray}\label{s2.2}
a(\varphi,v)= <f,v>,~~~\forall v\in H^{1}(\Omega).
\end{eqnarray}
\indent The Neumann eigenvalue problem associated with $n(x)$ is to find $k^{2}\in\mathbb{C}$ and a nontrivial $u$ such that
\begin{eqnarray}\label{s2.3}
\Delta u +k^{2} n(x)u=0~~in~ \Omega,~~~~
\frac{\partial u}{\partial \gamma}=0~~ on~ \partial\Omega.
\end{eqnarray}
\indent In this paper, we always assume $k^2$ is not an interior Neumann eigenvalue of (\ref{s2.3}). Under this assumption,  according to \cite{cakoni1}
 the Neumann-to-Dirichlet map $T: L^{2}(\partial\Omega)\to L^{2}(\partial\Omega)$ can be defined as follows. Let $f\in L^{2}(\partial\Omega)$, define $A: L^{2}(\partial\Omega)\to H^{1}(\Omega)$ by
\begin{eqnarray}\label{s2.4}
a(Af,v)= <f,v>,~~~\forall v\in H^{1}(\Omega),
\end{eqnarray}
and
$
Tf=(Af)'
$,
where $'$ denotes the restriction to $\partial\Omega$.
And (\ref{s2.1}) can be stated as the operator form:
\begin{eqnarray}\label{s2.5}
Tw=\mu w.
\end{eqnarray}
(\ref{s2.1}) and (\ref{s2.5}) are equivalent, namely, if $(\mu,w)\in \mathbb{C}\times L^{2}(\partial\Omega)$ is an eigenpair of (\ref{s2.5}), then $(\lambda,Aw)$ is an eigenpair of (\ref{s2.1}),
$\lambda=-\mu^{-1}$; conversely, if $(\lambda,u)$ is an
eigenpair of (\ref{s2.1}), then $(\mu,u')$
is an eigenpair of (\ref{s2.5}), $\mu=-\lambda^{-1}$.\\
\indent From \cite{cakoni1} we know
$T: L^{2}(\partial\Omega)\to L^{2}(\partial\Omega)$ is compact.
If $n(x)$ is real, then $T$ is also self-adjoint.\\
\indent Consider the dual problem of (\ref{s2.1}): Find
$\lambda^{*}\in \mathbb{C}$, $u^{*}\in
H^{1}(\Omega)\setminus \{0\}$ such that
\begin{eqnarray}\label{s2.6}
a(v, u^{*}) =-\overline{\lambda^{*}} <v, u^{*}>,~~~\forall v\in H^{1}(\Omega).
\end{eqnarray}
\indent The source problem associated with (\ref{s2.6}) is as
follows: Find $\varphi^*\in H^{1}(\Omega)$ such that
\begin{eqnarray}\label{s2.7}
a(v,\varphi^*)= <v,g>,~~~\forall v\in H^{1}(\Omega),
\end{eqnarray}
Define the corresponding Neumann-to-Dirichlet operator operator $T^{*}: L^{2}(\partial\Omega)$
$\to L^{2}(\partial\Omega)$ by
\begin{eqnarray}\label{s2.8}
&&a(v, A^{*}g) =<v,g>,~~~ \forall v\in
H^{1}(\Omega),
\end{eqnarray}
and $T^{*}g=( A^{*}g)'$. Then (\ref{s2.6}) has the equivalent operator form:
\begin{eqnarray}\label{s2.9}
T^{*}u^{*}=-\lambda^{*-1}u^{*}.
\end{eqnarray}
\indent It can be proved that $T^{*}$ is the adjoint operator of $T$
in the sense of inner product $<\cdot,\cdot>$. In fact, from
(\ref{s2.4}) and (\ref{s2.8}) we have
\begin{eqnarray*}
<Tf,g>=a(Af,A^*g)=<f,A^*g>=<f,T^*g>,~~~\forall
f,g\in L^{2}(\partial\Omega
).
\end{eqnarray*}
\indent Note that since $T^{*}$ is the adjoint operator of $T$, the
primal and dual eigenvalues are connected via
$\lambda=\overline{\lambda^{*}}$.\\
Let $\pi_{h}=\{\kappa\}$ be a regular $d$-simplex partition of $\Omega$
 (see
\cite{ciarlet}, pp. 131). We denote $h=\max_{\kappa\in \pi_{h}}h_{\kappa}$
where $h_{\kappa}$ is the diameter of element $\kappa$.
Let $\mathcal{E}_{h}$ denote the set of all $(d-1)$-faces of elements $\kappa\in\pi_{h}$. We split this set as follows: $\mathcal{E}_{h}
 =\mathcal{E}_{h}^{i}\cup\mathcal{E}_{h}^{b}$,
 with $\mathcal{E}_{h}^{i}$ and $\mathcal{E}_{h}^{b}$ being the sets of inner and boundary edges, respectively.
Let $S^{h}$
be the C-R element space defined on $\pi_{h}$:\\
\indent~~~ $S^{h}=\{v\in L^{2}(\Omega):v\mid_{\kappa} \in P_{1}(\kappa)$, $v$ is
continuous at the barycenters\\
\indent\indent~~~ ~~~~~~ of the $(d-1)$-faces of element $\kappa,~\forall \kappa\in\pi_{h}\}$.\\
The C-R element approximation of
(\ref{s2.1}) is: Find $\lambda_{h}\in \mathbb{C}$, $u_{h} \in
S^{h}\setminus\{0\}$, such that
\begin{eqnarray}\label{s2.10}
a_{h}(u_{h},v)=-\lambda_{h} <u_{h},v>,~~~\forall v\in S^{h},
\end{eqnarray}
where
$
a_{h}(u_{h},v)=\sum\limits_{\kappa\in\pi_{h}}\int\limits_{\kappa}(\nabla
u_{h}\cdot \nabla\bar{v}-k^{2}n(x) u_{h}\bar{v}) dx
$.\\
\indent Define $\|v\|_{h}=(\sum\limits_{\kappa\in
\pi_{h}}\|v\|_{1,\kappa}^{2})^{\frac{1}{2}}$,
$\|v\|_{1,\kappa}^{2}=\int_{\kappa}(\sum\limits_{i=1}^{d}|\frac{\partial v }{\partial
x_{i}}|^{2}+|v|^{2})dx$.
Evidently, $\|\cdot\|_{h}$ is the norm on $S^{h}$ and it is easy to know that $a_{h}(\cdot,\cdot)$ is not uniformly $S^{h}$-elliptic.\\
The C-R element approximation of (\ref{s2.2})
is: Find $\varphi_{h} \in S^{h}$, such that
\begin{eqnarray}\label{s2.11}
a_{h}(\varphi_{h},v)=<f,v>,~~~\forall v \in S^{h}.
\end{eqnarray}

\indent Since $k^2$ is not an interior Neumann eigenvalue of (\ref{s2.3}), from spectral approximation theory \cite{chatelin} we know that when $h$ is properly small  $k^2$ is also not a C-R element eigenvalue for (\ref{s2.3}).
So the discrete source problem (\ref{s2.11}) is uniquely solvable.
 Thus, we can define the discrete
operator $A_{h}:L_{2}(\partial\Omega)\to S^{h}$, satisfying
\begin{eqnarray}\label{s2.12}
a_{h}(A_{h}f,v)=<f,v>,~~~ \forall v \in S^{h}.
\end{eqnarray}
Let us denote by $\delta S^{h}$ the function space defined on
$\partial\Omega$, which are restriction of functions in $S^{h}$ to
$\partial\Omega$.
Define the discrete operator $T_{h}: L_{2}(\partial\Omega)\to \delta S^{h}\subset L_{2}(\partial\Omega)$, satisfying
$
T_{h}f=(A_{h}f)'
$.
Then (\ref{s2.10}) has the equivalent operator form:
\begin{eqnarray}\label{s2.13}
T_{h}w_{h}=\mu_{h} w_{h},
\end{eqnarray}
namely, if $(\mu_{h},w_{h})\in \mathbb{C}\times L^{2}(\partial\Omega)$ is an eigenpair of (\ref{s2.13}), then $(\lambda_{h},A_{h}w_{h})$ is an eigenpair of (\ref{s2.10}),
$\lambda_{h}=-\mu_{h}^{-1}$; conversely, if $(\lambda_{h},u_{h})$ is an
eigenpair of (\ref{s2.10}), then $(\mu_{h},u_{h}')$
is an eigenpair of (\ref{s2.13}), $\mu_{h}=-\lambda_{h}^{-1}$.\\
\indent The non-conforming finite element approximation of
(\ref{s2.6}) is given by:
 Find
$\lambda_{h}^{*} \in \mathbb{C}$, $u_{h}^{*}\in
S^{h}\setminus \{0\}$ such that
\begin{eqnarray}\label{s2.14}
a_{h}(v, u_{h}^{*}) =-\overline{\lambda_{h}^{*}}
<v,u_{h}^{*}>,~~~\forall v\in
S^{h}.
\end{eqnarray}
The C-R element approximation of (\ref{s2.7})
is: Find $\varphi_{h}^{*} \in S^{h}$, such that
\begin{eqnarray}\label{s2.15}
a_{h}(v,\varphi_{h}^{*})=<v,g>,~~~\forall v \in S^{h}.
\end{eqnarray}
Define the discrete operator $A_{h}^{*}: L^{2}(\partial\Omega)\to
S^{h}$ satisfying
\begin{eqnarray}\label{s2.16}
&&a_{h}(v, A_{h}^{*}g)=<v,g>,~~~\forall~v\in
S^{h},
\end{eqnarray}
and denote $T_{h}^{*}g=(A_{h}^{*}g)'$, then (\ref{s2.14}) has the following equivalent operator form
\begin{eqnarray}\label{s2.17}
T_{h}^{*}u_{h}^{*}=-\lambda_{h}^{*-1}u_{h}^{*}.
\end{eqnarray}
\indent It can be proved that $T_{h}^{*}$ is the adjoint operator of
$T_{h}$ in the sense of inner product $<\cdot,\cdot>$.
Hence, the primal and dual eigenvalues are connected via
$\lambda_{h}=\overline{\lambda_{h}^{*}}$.\\
\indent We need the following regularity estimates which play an important role in our theoretical analysis.
Note that for $v\in  H^{\frac{1}{2}}(\partial\Omega)$, $<f,v>$ has a continuous extension, still denoted by $<f,v>$, to $f\in  H^{-\frac{1}{2}}(\partial\Omega)$.\\

\noindent{\bf Lemma 1~} For any $f\in  H^{-\frac{1}{2}}(\partial\Omega)$,
let $<f,v>$ be
the dual product on $H^{-\frac{1}{2}}(\partial\Omega)\times H^{\frac{1}{2}}(\partial\Omega)$ in (\ref{s2.2}), then
 there exists a unique solution $\varphi\in H^{1}(\Omega)$ to (\ref{s2.2}) such that
\begin{eqnarray}\label{s2.18}
\|\varphi\|_{1}\lesssim \|f\|_{-\frac{1}{2},\partial\Omega}.
\end{eqnarray}
\noindent{\bf Proof~}
Since $k^2$ is not an interior Neumann eigenvalue of (\ref{s2.3}),
 there exists a unique solution $\varphi\in H^{1}(\Omega)$ to (\ref{s2.2}).
Denote
\begin{eqnarray*}
b(u,v)=\int_{\Omega}\nabla u\cdot\overline{\nabla v} +n(x)u\overline{v}dx.
\end{eqnarray*}
Referring to the proof of (14.11) in \cite{ciarlet}, we will show that,  for any $v\in H^1(\Omega)$, there exists a constant $C(\Omega)$ such that
\begin{eqnarray}\label{s2.19}
\|v\|_{1}\leq C(\Omega)\sqrt{Re b(v,v)}\equiv C(\Omega)(\int\limits_{\Omega}|\nabla v|^{2} +n_{1}(x)|v|^{2}dx)^{\frac{1}{2}}.
\end{eqnarray}
From (\ref{s2.19}), it is easy to verify that $\sqrt{Re b(v,v)}$ is a norm on $H^{1}(\Omega)$ that is equivalent to the norm $\|\cdot\|_{1}$.\\
If (\ref{s2.19}) is false, there exists a sequence $\{v_{i}\}\subset H^{1}(\Omega)$ such that
 \begin{eqnarray}\label{s2.20}
\|v_{i}\|_{1}=1~~for~all~i\geq 1,~~~\lim\limits_{i\to\infty}(\int\limits_{\Omega}|\nabla v_i|^{2} +n_{1}(x)|v_i|^{2}dx)^{\frac{1}{2}}=0.
\end{eqnarray}
Since the $\{v_{i}\}$ is bounded in $H^{1}(\Omega)$, there exists a subsequence, again denoted $\{v_{i}\}$ for notational convenience,
that converges in $L^{2}(\Omega)$. Since
$\lim\limits_{i\to\infty}|v_{i}|_{1}=0$, by (\ref{s2.20}), and since $H^{1}(\Omega)$
is complete, the sequence $\{v_{i}\}$ converges in $H^{1}(\Omega)$. The limit $v$ of this sequence is such that
\begin{eqnarray*}
|v|_{1}=\lim\limits_{i\to\infty}|v_{i}|=0,
\end{eqnarray*}
and thus, $v$ is a constant. Using (\ref{s2.20}), we have
 \begin{eqnarray*}
(\int\limits_{\Omega}n_{1}(x)|v|^{2}dx)^{\frac{1}{2}}=\lim\limits_{i\to\infty}(\int\limits_{\Omega}n_{1}(x)|v_{i}|^{2}dx)^{\frac{1}{2}}=0,
\end{eqnarray*}
hence we conclude that $n_{1}(x)=0$. But this contradicts the hypothesis $n_{1}(x)>0$. It implies that (\ref{s2.19}) holds.\\
(\ref{s2.3}) can be rewritten as: Find $\tilde{\lambda}\in
\mathbb{C}$, $u \in H^{1}(\Omega)\backslash\{0\}$ such that
\begin{eqnarray}\label{s2.21}
b(u,v)=\tilde{\lambda} (nu,v),~~~\forall v\in H^{1}(\Omega),
\end{eqnarray}
Since $k^2$ is not an interior Neumann eigenvalue for (\ref{s2.3}), $k^2+1$ is not an eigenvalue
of (\ref{s2.21}).
Define the map $B: H^{1}(\Omega)\to H^{1}(\Omega)$ by
\begin{eqnarray}\label{s2.22}
b(Bg,v)= (ng,v),~~~\forall v\in H^{1}(\Omega).
\end{eqnarray}
Then (\ref{s2.21}) has the operator form:
\begin{eqnarray*}
Bu=\tilde{\lambda}^{-1} u.
\end{eqnarray*}
And $B$ is compact, $\frac{1}{k^2+1}$ is not an eigenvalue
of $B$.
So $(B-\frac{1}{k^2+1}I)^{-1}:H^{1}(\Omega)\to H^{1}(\Omega)$ is bounded.
Let $\psi \in H^{1}(\Omega)$ be the solution of the following equation:
 \begin{eqnarray}\label{s2.23}
b(\psi,v)=<f,v>,~~~\forall v\in H^{1}(\Omega),
\end{eqnarray}
then we have
$
\|\psi\|_{1}\lesssim\|f\|_{-\frac{1}{2},\partial\Omega}
$.
From (\ref{s2.22}) we obtain
 \begin{eqnarray*}
 &&a(\varphi,v)=b(\varphi,v)-(k^2+1)(n\varphi,v)\\
 &&~~~=b(\varphi,v)-(k^2+1)b(B\varphi,v)=b((I-(k^2+1)B)\varphi,v),
  \end{eqnarray*}
which together with  (\ref{s2.2}) and (\ref{s2.23}) yields
 \begin{eqnarray*}
 b((I-(k^2+1)B)\varphi,v)=b(\psi,v),~~~\forall v\in H^{1}(\Omega).
 \end{eqnarray*}
 Thus we have
  \begin{eqnarray*}
 &&\|\varphi\|_{1}=\|(\frac{1}{k^2+1}I-B)^{-1}\frac{1}{k^2+1}\psi\|_{1}\\
 &&~~~\leq\|(\frac{1}{k^2+1}I-B)^{-1}\|_{1}|\frac{1}{k^2+1}|\|\psi\|_{1}
 \lesssim\|f\|_{-\frac{1}{2},\partial\Omega},
 \end{eqnarray*}
 and the proof is completed.~~~$\Box$\\

\noindent{\bf Lemma 2~}Assume that $\Omega\subset \mathbb{R}^{2}$ is a polygonal with $\omega$ being the largest interior angle, and $\varphi$ is the solution of (\ref{s2.2}).\\
i)~If $f\in L^{2}(\partial\Omega)$, then $\varphi\in H^{1+\frac{r}{2}}(\Omega)$ and
\begin{eqnarray}\label{s2.24}
\|\varphi\|_{1+\frac{r}{2}}\leq C_{\Omega}\|f\|_{0,\partial\Omega}.
\end{eqnarray}
ii)~If $f\in H^{\frac{1}{2}}(\partial\Omega)$, then
$\varphi\in
H^{1+r}(\Omega)$ satisfying
\begin{eqnarray}\label{s2.25}
 \|\varphi\|_{1+r}\leq
C_{\Omega}\|f\|_{\frac{1}{2},\partial\Omega}.
\end{eqnarray}
Here $r=1$ when $\omega<\pi$, and $r<\frac{\pi}{\omega}$ when $\omega>\pi$,
 and $C_{\Omega}$ is a priori constant.\\
\noindent{\bf Proof~}
 Consider the auxiliary boundary value problem:
\begin{eqnarray}\label{s2.26}
&&\Delta \varphi_{1}+\varphi_{1}=0,~~~\frac{\partial\varphi_{1}}{\partial\gamma}=f,\\\label{s2.27}
&&\Delta \varphi_{2}+\varphi_{2}=-k^2n(x)(\varphi_{1}+\varphi_{2})+\varphi_{1}+\varphi_{2},~~~\frac{\partial\varphi_{2}}{\partial\gamma}=0.
\end{eqnarray}
Let $\varphi_{1}$ and $\varphi_{2}$ be the solution of (\ref{s2.26}) and (\ref{s2.27}), respectively, then it is easy to see that $\varphi=\varphi_{1}+\varphi_{2}$.
Since $\Omega\subset\mathbb{R}^{2}$,
from classical regularity results (see \cite{dauge}, or Proposition 4.1 in \cite{alonso} and Proposition 4.4 in \cite{bermudez}) we have
\begin{eqnarray*}
\|\varphi_{1}\|_{1+(\frac{1}{2}+s)r}\lesssim \|f\|_{s,\partial\Omega},~~~s=0,\frac{1}{2},
\end{eqnarray*}
and from classical regularity result for the Laplace problem with homogeneous Neumann boundary condition we have
\begin{eqnarray*}
\|\varphi_{2}\|_{1+r}\lesssim \|-k^2n(x)(\varphi_{1}+\varphi_{2})+\varphi_{1}+\varphi_{2}\|_{0,\Omega}.
\end{eqnarray*}
Thus we get
\begin{eqnarray*}
&&\|\varphi\|_{1+(\frac{1}{2}+s)r}\lesssim \|\varphi_{1}\|_{1+(\frac{1}{2}+s)r}+\|\varphi_{2}\|_{1+(\frac{1}{2}+s)r}
\lesssim \|f\|_{s,\partial\Omega}+\|\varphi\|_{0},~~~s=0,\frac{1}{2}.
\end{eqnarray*}
Substituting (\ref{s2.18}) into the above inequality we get (\ref{s2.24}) and (\ref{s2.25}).~~~$\Box$\\

\noindent{\bf Remark 1}(Regularity in $\mathbb{R}^3$)~
When $\Omega\subset \mathbb{R}^{3}$ is a polyhedral domain, regularity of the solution of the Neumann problem (\ref{s2.26}) has been discussed by many scholars.
Referring Theorem 4 in \cite{savare}, Remark 2.1 in \cite{garau}, \cite{jerison} and \cite{dauge},
and using the argument of Lemma 2 in this paper,
we think the following regularity hypothesis $R(\Omega)$ is reasonable:\\
\indent{\bf Hypothesis $ R(\Omega)$~}
  Let $\varphi$ be the solution of (\ref{s2.2}) with
$f\in L^{2}(\partial\Omega)$. When
$\Omega\subset \mathbb{R}^{3}$, we have $\varphi\in H^{1+r3}(\Omega)$ for all
 $r3\in (0,\frac{1}{2})$ and
\begin{eqnarray}\label{s2.28}
\|\varphi\|_{1+r3}\leq C_{\Omega}\|f\|_{0,\partial\Omega}.
\end{eqnarray}

It is easy to know that Lemmas 1-2 and Remark 1 are also valid for the dual problem (\ref{s2.7}).

\section{The consistency term and the extension of Strang lemma}

\indent Define $S^{h} +H^{1}(\Omega)=\{ w_h+
w: w_h\in S^h,  w\in H^{1}(\Omega)\}.$\\
Let
$\varphi$ and $\varphi^{*}$ be the solution of
(\ref{s2.2}) and (\ref{s2.7}), respectively. Define the
consistency terms: For any $v\in S^{h} +H^{1}(\Omega)$,
\begin{eqnarray}\label{s3.1}
&&D_{h}(\varphi,v)=a_{h}(\varphi,v)-<f, v>,\\\label{s3.2}
&&D_{h}^{*}(v,\varphi^{*})=a_{h}(v,\varphi^{*})-<v,g>.
\end{eqnarray}
\indent In order to analyze error estimates of the consistency terms, we need the following trace inequalities.\\

\noindent{\bf Lemma 3~}
For any $\kappa\in \pi_{h}$, the following trace inequalities hold:
\begin{eqnarray}
&&\|w\|_{0,\partial\kappa}
\lesssim h_{\kappa}^{-\frac{1}{2}}\|w\|_{0,\kappa}+h_{\kappa}^{\frac{1}{2}}|w|_{1,\kappa},\nonumber\\\label{s3.3}
&&\|\nabla w\|_{0,\partial\kappa}
\lesssim h_{\kappa}^{-\frac{3}{2}}\| w\|_{0,\kappa}+h_{\kappa}^{-\frac{1}{2}}| w|_{1,\kappa}+h_{\kappa}^{s-\frac{1}{2}}|w|_{1+s,\kappa}~(\frac{1}{2}\leq s\leq 1).
\end{eqnarray}
\noindent{\bf Proof~}
 The conclusion is followed by using the trace theorem on the reference element and the scaling argument (see, e.g., Lemma 2.2 in \cite{yang1}).~~~$\Box$\\

\indent The following Green's formula (see (2.7) in \cite{caiz}, (3.11) in \cite{bernardi}
and Corollary 2.2 in \cite{girault}) will play a crucial role in our analysis:
\begin{eqnarray}\label{s3.4}
\int\limits_{\partial\kappa}(\nabla w\cdot\gamma)vds
=\int\limits_{\kappa}\Delta w vdx+\int\limits_{\kappa}\nabla w\cdot\nabla vdx~~~\forall \kappa\in\pi_{h},
\end{eqnarray}
where $w\in H^{1+\epsilon}(\kappa)$ with $\Delta w\in L^{2}(\kappa)$ and $v\in H^{1-\epsilon}(\kappa)$ with
$0\leq\epsilon<\frac{1}{2}$.\\

\noindent{\bf Lemma 4~}Suppose that
$\varphi\in H^{1+t}(\Omega) (0<t<\frac{1}{2})$ is the solution of (\ref{s2.2}) and {\bf Hypothesis $R(\Omega)$} holds,
then
\begin{eqnarray}\label{s3.5}
\|\nabla \varphi\cdot\gamma\|_{t-\frac{1}{2},\ell}
\lesssim\|\varphi\|_{1+t,\kappa}~~~\forall \kappa\in
\pi_{h},~\ell\in\partial\kappa.
\end{eqnarray}
\noindent{\bf Proof~}
Inequality (\ref{s3.5}) is contained in the proof of Corollary 3.3 on page 1384
of \cite{bernardi}, see also Lemma 2.1 in \cite{caiz}. For the convenience of readers, we write the proof here.\\
For any $g\in H^{\frac{1}{2}-t}(\ell)$,
it is proven by going to a reference
element and using the inverse trace theorem (see page 387 in \cite{kufner}, page 1767 in \cite{caiz})) that
there exists a lifting $w_{g}$ of $g$ such that $w_{g}\in H^{1-t}(\kappa)$, $w_{g}|_{\ell}=g$, $w_{g}|_{\partial\kappa\setminus\ell}=0$,
 and
\begin{eqnarray}\label{s3.6}
\|\nabla w_{g}\|_{-t,\kappa}+h_{\kappa}^{t-1}\| w_{g}\|_{0,\kappa}\leq C\|g\|_{\frac{1}{2}-t,\ell}.
\end{eqnarray}
From Green's formula (\ref{s3.4}), (\ref{s2.2}), Cauchy-Schwarz inequality,
the definition of the dual norm and (\ref{s3.6}) we deduce
\begin{eqnarray*}
&&\int\limits_{\ell}\nabla\varphi\cdot\gamma\bar{g}ds=\int\limits_{\ell}\nabla\varphi\cdot\gamma\overline{w_{g}}ds
=\int\limits_{\kappa}\Delta\varphi\overline{w_{g}}dx+\int\limits_{\kappa}\nabla\varphi\cdot\overline{\nabla w_{g}}dx\nonumber\\
&&~~~=\int\limits_{\kappa}-k^{2}n\varphi\overline{w_{g}}dx+\int\limits_{\kappa}\nabla\varphi\cdot\overline{\nabla w_{g}}dx\lesssim\|-k^{2}n\varphi\|_{0,\kappa}\|\overline{w_{g}}\|_{0,\kappa}+\|\nabla\varphi\|_{t,\kappa}\|\overline{\nabla w_{g}}\|_{-t,\kappa}\nonumber\\
&&~~~\lesssim\|-k^{2}n\varphi\|_{0,\kappa}h_{\kappa}^{1-t}\|g\|_{\frac{1}{2}-t,\ell}+\|\nabla\varphi\|_{t,\kappa}\|g\|_{\frac{1}{2}-t,\ell}
\lesssim \|\varphi\|_{1+t,\kappa}\|g\|_{\frac{1}{2}-t,\ell},
\end{eqnarray*}
thus by the definition of the dual norm we obtain
\begin{eqnarray*}
\|\nabla \varphi\cdot\gamma\|_{t-\frac{1}{2},\ell}=\sup\limits_{g\in H^{\frac{1}{2}-t}(\ell)}\frac{|\int\limits_{\ell}\nabla\varphi\cdot\gamma\bar{g}ds|}
{\|g\|_{\frac{1}{2}-t,\ell}}\lesssim\|\varphi\|_{1+t,\kappa}.
\end{eqnarray*}
This completes the proof of the lemma.~~~$\Box$\\

Based on the standard argument (see, e.g., \cite{alonso,liq2,yang1}), the following consistency error estimates will be
proved.\\

\noindent{\bf Theorem 1~} Let $\varphi$ and $\varphi^{*}$ be the solution of (\ref{s2.2}) and (\ref{s2.7}), respectively, and suppose that $\varphi, \varphi^{*}\in H^{1+t}(\Omega)$ with $t\in [s,1]$ and {\bf Hypothesis $R(\Omega)$} holds, then
\begin{eqnarray}\label{s3.7}
&&|D_{h}(\varphi,v)|\lesssim h^{t}\|\varphi\|_{1+t}\|v\|_{h},~~~\forall v\in
S^{h}+H^{1}(\Omega),\\\label{s3.8}
&&|D_{h}^*(v,\varphi^*)|\lesssim h^{t}\|\varphi^*\|_{1+t}\|v\|_{h},~~~\forall v\in
S^{h}+H^{1}(\Omega),
\end{eqnarray}
where s=$\frac{r}{2}$ when $\Omega\subset\mathbb{R}^2$, s=$r3$ when $\Omega\subset\mathbb{R}^3$.\\
\noindent{\bf Proof~}
 Let
$[[\cdot]]$
denote the jump across an inner face $\ell\in \mathcal{E}_{h}^{i}$. By
Green's formula (\ref{s3.4}) we deduce
\begin{eqnarray}\label{s3.9}
&&D_{h}(\varphi,v)=a_{h}(\varphi,v)-<f,v>=\sum\limits_{\kappa\in\pi_{h}}\int\limits_{\kappa}(-\Delta \varphi-k^{2}n(x)\varphi)\bar{v}dx\nonumber\\
&&~~~~~~+\sum\limits_{\kappa\in\pi_{h}}\int\limits_{\partial\kappa}\frac{\partial\varphi}{\partial\gamma}\bar{v}ds
-\int\limits_{\partial\Omega}\frac{\partial\varphi}{\partial\gamma}\bar{v}ds=\sum\limits_{\ell\in \mathcal{E}_{h}^{i}}
\int_{\ell}\frac{\partial \varphi}{\partial \gamma}\bar{[[v]]}ds.
\end{eqnarray}
Let $\ell$ be a ($d-1$)-face of $\kappa$, define
\begin{eqnarray*}
P_{\ell}f=\frac{1}{|\ell|}\int\limits_{\ell}fds,
~~~P_{\kappa}f=\frac{1}{|\kappa|}\int\limits_{\kappa}fdx.\\
\end{eqnarray*}
For $\ell\in
\mathcal{E}_{h}^{i}$, suppose that $\kappa_{1},\kappa_{2}\in \pi_{h}$ such that
$\kappa_{1}\bigcap \kappa_{2}=\ell$. Since $[[v]]$ is a linear function vanishing
at the barycenters of $\ell$, we have
\begin{eqnarray}\label{s3.10}
&&|\int_{\ell}\frac{\partial \varphi}{\partial \gamma}\overline{[[v]]}ds|
=|\int_{\ell}(\frac{\partial \varphi}{\partial
\gamma}-P_{\ell}(\frac{\partial \varphi}{\partial \gamma}))\overline{[[v]]}ds|\nonumber\\
&&~~~=|\int_{\ell}(\frac{\partial
\varphi}{\partial \gamma}-P_{\ell}(\frac{\partial \varphi}{\partial
\gamma}))\overline{([[v]]-P_{\ell}[[v]])}ds|\nonumber\\
&&~~~= |\int_{\ell}\frac{\partial
\varphi}{\partial \gamma}
\overline{([[v]]-P_{\ell}[[v]])}ds|.
\end{eqnarray}
Then, when $t\in [ \frac{1}{2},1]$, using Schwarz inequality we deduce
\begin{eqnarray}\label{s3.11}
&&|\int_{\ell}\frac{\partial \varphi}{\partial \gamma}\overline{[[v]]}ds|
\leq \sum\limits_{i=1,2}\|\nabla \varphi\cdot \gamma-
P_{\ell}(\nabla \varphi\cdot
\gamma)\|_{0,\ell}\|v|_{\kappa_{i}}-P_{\ell}(v|_{\kappa_{i}})\|_{0,\ell}\nonumber\\
&&~~~ \leq\sum\limits_{i=1,2}
\|\nabla (\varphi-\varphi_{I})\cdot \gamma\|_{0,\ell}\|v|_{\kappa_{i}}-P_{\kappa_{i}}(v|_{\kappa_{i}})\|_{0,\ell},
\end{eqnarray}
by Lemma 3 and the standard error estimates for
$L^{2}$-projection,  we deduce
\begin{eqnarray*}
\|\nabla (\varphi-\varphi_{I})\cdot \gamma\|_{0,\ell}\lesssim h^{t-\frac{1}{2}}||\varphi||_{1+t,\kappa_{i}},~~~
\|v|_{\kappa_{i}}-P_{\kappa_{i}}(v|_{\kappa_{i}})\|_{0,\ell}\lesssim h^{\frac{1}{2}}||v||_{1,\kappa_{i}}.
\end{eqnarray*}
Substituting the above two estimates into (\ref{s3.11}), we obtain
\begin{eqnarray}\label{s3.12}
|\int_{\ell}\frac{\partial \varphi}{\partial \gamma}\overline{[[v]]}ds| \lesssim\sum\limits_{i=1,2}h^{t}\|\varphi\|_{1+t,\kappa_{i}}\|v\|_{1,\kappa_{i}},
\end{eqnarray}
and substituting (\ref{s3.12}) into (\ref{s3.9}) we conclude that
(\ref{s3.7}) holds.\\
When $t<\frac{1}{2}$, from (\ref{s3.10}) we deduce that
\begin{eqnarray}\label{s3.13}
|\int_{\ell}\frac{\partial \varphi}{\partial \gamma}\overline{[[v]]}ds|\leq
\|\nabla \varphi\cdot \gamma\|_{t-\frac{1}{2},\ell}\|[[v]]-P_{\ell}[[v]]\|_{\frac{1}{2}-t,\ell}.
\end{eqnarray}
By using inverse estimate, Lemma 3 and the error estimate of $L^{2}$-projection, we derive
\begin{eqnarray*}
&&\|[[v]]-P_{\ell}[[v]]\|_{\frac{1}{2}-t,\ell}\lesssim h_{\ell}^{t-\frac{1}{2}}\|[[v]]-P_{\ell}[[v]]\|_{0,\ell}
\lesssim  \sum\limits_{i=1,2}h_{\kappa_{i}}^{t}||v||_{1,\kappa_{i}}.
\end{eqnarray*}
Substituting the above estimate and (\ref{s3.5}) into (\ref{s3.13}), we obtain
\begin{eqnarray*}
|\int_{\ell}\frac{\partial \varphi}{\partial \gamma}\overline{[[v]]}ds|
 \lesssim\sum\limits_{i=1,2}\|\varphi\|_{1+t,\kappa_{i}}h_{\kappa_{i}}^{t}||v||_{1,\kappa_{i}},
\end{eqnarray*}
plugging the above inequality into (\ref{s3.9}) we also get (\ref{s3.7}). \\
\indent Using the same argument as above,  we can prove
(\ref{s3.8}).~~~$\Box$\\

The C-R element approximation (\ref{s2.11}) of (\ref{s2.2}) does not satisfy the condition of Strang lemma, that is $a_{h}(\cdot,\cdot)$
is not uniformly $S^{h}$-elliptic.
To overcome this difficulty,
 Inspired by the work in \S 5.7 in \cite{brenner}, next we use the standard duality technique to prove an
extension version of the well-known Strang lemma. \\
First, we will use the standard duality argument to prove that $\|\varphi-\varphi_{h}\|_{0}$ is a quantity of higher order than $\|\varphi-\varphi_{h}\|_{h}$.\\
Introduce the auxiliary problem: Find $\psi\in H^{1}(\Omega)$, such that
\begin{eqnarray}\label{s3.14}
a(v,\psi)=(v,g),~~~~\forall~ v\in H^{1}(\Omega).
\end{eqnarray}
Let $\psi$ be the solution of (\ref{s3.14}), then
from elliptic regularity estimates for homogeneous Neumann boundary value problem we know that there exists $r_{N}>0$, such that
\begin{eqnarray}\label{s3.15}
\|\psi\|_{1+r_{N}}\lesssim \|g\|_{0}.
\end{eqnarray}
Let $\hat{E}_{h}(v,\psi)=a_{h}(v,\psi)-(v,g)$,  then
\begin{eqnarray}\label{s3.16}
|\hat{E}_{h}(v,\psi)|\lesssim h^{r_{N}}\|\psi\|_{1+r_{N}}\|v\|_{h},~~~\forall v\in H^{1}(\Omega)+S^{h}.
\end{eqnarray}

\noindent{\bf Lemma 5~}
 Let $\varphi$ and $\varphi_{h}$ be the solution of
(\ref{s2.2}) and
(\ref{s2.11}), respectively, and let $\varphi^{*}$
and $\varphi_{h}^{*}$ be the solution of
 (\ref{s2.7}) and (\ref{s2.15}), respectively,
then
\begin{eqnarray}\label{s3.17}
&&\|\varphi-\varphi_{h}\|_{0}\lesssim h^{r_{N}}\|\varphi-\varphi_{h}\|_{h},\\\label{s3.18}
&&\|\varphi^{*}-\varphi_{h}^{*}\|_{0}\lesssim  h^{r_{N}}\|\varphi^{*}-\varphi_{h}^{*}\|_{h}.
\end{eqnarray}
\noindent{\bf Proof~}
By Riesz representation theorem we have
\begin{eqnarray}\label{s3.19}
\|\varphi-\varphi_{h}\|_{0}=\sup\limits_{g\in L^2(\Omega), g\not=0}\frac{|(\varphi-\varphi_{h},g)|}{\|g\|_{0}}.
\end{eqnarray}
Let $\psi_{I}\in S^{h}$ be the C-R element interpolation function of $\psi$, the solution of (3.13),
 then according to the interpolation theory (see \cite{ciarlet}) we have
\begin{eqnarray}\label{s3.20}
\|\psi-\psi_{I}\|_{h}\lesssim  h^{r_{N}}\|\psi\|_{1+r_{N}}.
\end{eqnarray}
By computing, we deduce
\begin{eqnarray*}
&&<f,\psi-\psi_{I}>=<f,\psi>-<f,\psi_{I}>=a(\varphi,\psi)-a_{h}(\varphi_{h},\psi_{I})\nonumber\\
&&~~~=a(\varphi,\psi)-a_{h}(\varphi,\psi_{I})+a_{h}(\varphi,\psi_{I})-a_{h}(\varphi_{h},\psi_{I})\nonumber\\
&&~~~=a_{h}(\varphi,\psi-\psi_{I})+a_{h}(\varphi-\varphi_{h},\psi_{I})\nonumber\\
&&~~~=a_{h}(\varphi,\psi-\psi_{I})+a_{h}(\varphi-\varphi_{h},\psi_{I}-\psi)+a_{h}(\varphi-\varphi_{h},\psi),
\end{eqnarray*}
and
\begin{eqnarray*}
&&a_{h}(\varphi,\psi-\psi_{I})+a_{h}(\varphi-\varphi_{h},\psi)=a_{h}(\varphi,\psi-\psi_{I})-<f,\psi-\psi_{I}>\\
&&~~~~~~+<f,\psi-\psi_{I}>+a_{h}(\varphi-\varphi_{h},\psi)-(\varphi-\varphi_{h},g)+(\varphi-\varphi_{h},g)\\
&&~~~=D_{h}(\varphi,\psi-\psi_{I})+\hat{E}_{h}(\varphi-\varphi_{h},\psi)+<f,\psi-\psi_{I}>+(
\varphi-\varphi_{h},g).
\end{eqnarray*}
Combining the above two inequalities we get
\begin{eqnarray*}
(\varphi-\varphi_{h},g)=a_{h}(\varphi-\varphi_{h},\psi-\psi_{I})-D_{h}(\varphi,\psi-\psi_{I})-\hat{E}_{h}(\varphi-\varphi_{h},\psi).
\end{eqnarray*}
Substituting the above equality into (\ref{s3.19}) we get
\begin{eqnarray}\label{s3.21}
&&\|\varphi-\varphi_{h}\|_{0}\leq
\sup\limits_{g\in L^2(\Omega), g\not=0}\frac{|a_{h}(\varphi-\varphi_{h},\psi-\psi_{I})|}{\|g\|_{0}}\nonumber\\
&&~~~~~~+\sup\limits_{g\in L^2(\Omega), g\not=0}\frac{|D_{h}(\varphi,\psi-\psi_{I})+\hat{E}_{h}(\varphi-\varphi_{h},\psi)|}{\|g\|_{0}}.
\end{eqnarray}
Let $I_{h}^{C}$ be the Lagrange interpolation operator, then $I_{h}^{C}\psi\in H^{1}(\Omega)\cap S^{h}$.
According to (\ref{s3.1}) and (\ref{s2.2})
we have $D_{h}(\varphi,v)=0$, $\forall v\in H^{1}(\Omega)$, and we deduce
\begin{eqnarray}\label{s3.22}
&&|D_{h}(\varphi,\psi-\psi_{I})|=|D_{h}(\varphi,I_{h}^{C}\psi-\psi_{I})|=|a_{h}(\varphi-\varphi_{h}, I_{h}^{C}\psi-\psi_{I} )|\nonumber\\
&&~~~\lesssim \|\varphi-\varphi_{h}\|_{h}\|I_{h}^{C}\psi-\psi_{I}\|_{h}\lesssim h^{r_{N}} \|\psi\|_{1+r_{N}}\|\varphi-\varphi_{h}\|_{h}.
\end{eqnarray}
Substituting (\ref{s3.16}), (\ref{s3.20}) and (\ref{s3.22}) into (\ref{s3.21})
 we obtain the desired result (\ref{s3.17}).\\
Using the same argument as (\ref{s3.17}) we can prove (\ref{s3.18}).~~~$\Box$\\

\indent Now we are ready to prove the following extension of Strang lemma.\\
\noindent{\bf Theorem 2~}
 Let $\varphi$ and $\varphi_{h}$ be the solution of
(\ref{s2.2}) and
(\ref{s2.11}), respectively, then
\begin{eqnarray}\label{s3.23}
&&\inf\limits_{v\in
S^{h}}\|\varphi-v\|_{h}+\sup\limits_{v\in
S^{h}\setminus\{0\}}\frac{
|D_{h}(\varphi,v)|}{\|v\|_{h}}\lesssim
\|\varphi-\varphi_{h}\|_{h}\nonumber\\
&&~~~\lesssim \inf\limits_{v\in
S^{h}}\|\varphi-v\|_{h}+\sup\limits_{v\in
S^{h}\setminus\{0\}} \frac{
|D_{h}(\varphi,v)|}{\|v\|_{h}}.
\end{eqnarray}
Let $\varphi^{*}$
and $\varphi_{h}^{*}$ be the solution of
 (\ref{s2.7}) and (\ref{s2.15}), respectively,
then
\begin{eqnarray}\label{s3.24}
&&\inf\limits_{v\in
S^{h}}\|\varphi^{*}-v\|_{h}+\sup\limits_{v\in
S^{h}\setminus\{0\}}\frac{
|D_{h}^*(v, \varphi^{*})|}{\|v\|_{h}}\lesssim
\|\varphi^{*}-\varphi_{h}^{*}\|_{h}\nonumber\\
&&~~~\lesssim \inf\limits_{v\in
S^{h}}\|\varphi^{*}-v\|_{h}+\sup\limits_{v\in
S^{h}\setminus\{0\}} \frac{| D_{h}^*(v,
\varphi^{*})|}{\|v\|_{h}}.
\end{eqnarray}
\noindent{\bf Proof~}
Denote
\begin{eqnarray*}
\mathcal{A}_{h}(u,v)=a_{h}(u,v)+K(u,v),~~~\forall u,v\in S^{h}+H^{1}(\Omega),
\end{eqnarray*}
where $K>\|k^{2}n\|_{0,\infty}$.
Then
we know that $\mathcal{A}_{h}$ satisfies the uniform $S^{h}$-ellipticity:
\begin{eqnarray*}
|\mathcal{A}_{h}(v,v)|\geq \min\{1, K-\|k^{2}n\|_{0,\infty}\}\|v\|_{h}^{2},~~~\forall v\in S^{h}.
\end{eqnarray*}
And thus, for any $v\in S^{h}$,
\begin{eqnarray*}
&&\|\varphi_{h}-v\|_{h}^{2}\lesssim| \mathcal{A}_{h}(\varphi_{h}-v, \varphi_{h}-v)|\\
&&~~~=C|a_{h}(\varphi-v,\varphi_{h}-v)+<f,\varphi_{h}-v>-a_{h}(\varphi,\varphi_{h}-v)+K\|\varphi_{h}-v\|_{0}^{2}|.
\end{eqnarray*}
When $\|\varphi_{h}-v\|_{h}\not= 0$, dividing both sides of the above by  $\|\varphi_{h}-v\|_{h}$ we obtain
\begin{eqnarray*}
&&\|\varphi_{h}-v\|_{h} \lesssim
\|\varphi-v\|_{h}+\frac{|a_{h}(\varphi, \varphi_{h}-v)
-<f,\varphi_{h}-v>|}{\|\varphi_{h}-v\|_{h}}+K\|\varphi_{h}-v\|_{0}\\
&&~~~\lesssim \|\varphi-v\|_{h}+\sup\limits_{v\in
S^{h}\setminus\{0\}}
\frac{|D_{h}(\varphi,v)|}{\|v\|_{h}}+\|\varphi-\varphi_{h}\|_{0}.
\end{eqnarray*}
From the triangular inequality and (\ref{s3.17}) we get
\begin{eqnarray*}
&&\|\varphi-\varphi_{h}\|_{h}\leq
\|\varphi-v\|_{h}+\|v-\varphi_{h}\|_{h}\nonumber\\
&&~~~\lesssim \|\varphi-v\|_{h}+\sup\limits_{v\in
S^{h}\setminus\{0\}}
\frac{|D_{h}(\varphi,v)|}{\|v\|_{h}}+h^{r_{N}}\|\varphi-\varphi_{h}\|_{h}.
\end{eqnarray*}
The second inequality of (\ref{s3.23}) is proved.\\
From
\begin{eqnarray*}
|a_{h}(\varphi-\varphi_{h}, v)|\leq
\|\varphi-\varphi_{h}\|_{h}\|v\|_{h},~~~~\forall~
v\in S^{h},
\end{eqnarray*}
we get
\begin{eqnarray*}
\|\varphi-\varphi_{h}\|_{h}\geq \frac{
|a_{h}(\varphi,v)-a_{h}(\varphi_{h},v)|}{\|v\|_{h}}=\frac{ |D_{h}(\varphi,v)|}{\|v\|_{h}},
\end{eqnarray*}
which together with
$\|\varphi-\varphi_{h}\|_{h}\geq
\inf\limits_{v\in S^{h}}\|\varphi-v\|_{h}$ we obtain
the first inequality of (\ref{s3.23}).\\
\indent Similarly we can prove (\ref{s3.24}). The proof is
completed.~~~$\Box$\\

\indent Now we can state the error estimates of C-R element approximation for (\ref{s2.2}) and (\ref{s2.7}).\\
\noindent{\bf Theorem 3~}
Under the conditions of Theorems 1 and 2,
we have
\begin{eqnarray}\label{s3.25}
&&\|\varphi-\varphi_{h}\|_{h}\leq
Ch^{t}\|\varphi\|_{1+t},\\\label{s3.26}
&&\|\varphi-\varphi_{h}\|_{0,\partial\Omega}\leq
Ch^{t+s}\|\varphi\|_{1+t}\\\label{s3.27}
&&\|\varphi^*-\varphi_{h}^*\|_{h}\leq
Ch^{t}\|\varphi^*\|_{1+t},\\\label{s3.28}
&&\|\varphi^*-\varphi_{h}^*\|_{0,\partial\Omega}\leq
Ch^{t+s}\|\varphi^*\|_{1+t}.
\end{eqnarray}
\noindent{\bf Proof~}
From  Theorem 2, the interpolation error estimate and Theorem 1 we can obtain (\ref{s3.25}) and (\ref{s3.27}).
From (\ref{s2.2}) and
(\ref{s2.11}) we deduce
\begin{eqnarray*}
&&<f,\varphi^{*}-\varphi_{h}^{*}>=<f,\varphi^{*}>-<f,\varphi_{h}^{*}>\nonumber\\
&&~~~=a(\varphi,\varphi^{*})-a_{h}(\varphi_{h},\varphi_{h}^{*})
=a(\varphi,\varphi^{*})-a_{h}(\varphi_{h},\varphi^{*})+a_{h}(\varphi_{h},\varphi^{*})-a_{h}(\varphi_{h},\varphi_{h}^{*})\nonumber\\
&&~~~=a_{h}(\varphi-\varphi_{h},\varphi^{*})+a_{h}(\varphi_{h},\varphi^{*}-\varphi_{h}^{*})\nonumber\\
&&~~~=a_{h}(\varphi-\varphi_{h},\varphi^{*})+a_{h}(\varphi,\varphi^{*}-\varphi_{h}^{*})-a_{h}(\varphi-\varphi_{h},\varphi^{*}-\varphi_{h}^{*}),
\end{eqnarray*}
and from (\ref{s3.1}) and (\ref{s3.2}) we get
\begin{eqnarray*}
&&a_{h}(\varphi-\varphi_{h},\varphi^{*})+a_{h}(\varphi,\varphi^{*}-\varphi_{h}^{*})
=a_{h}(\varphi-\varphi_{h},\varphi^{*})-<\varphi-\varphi_{h},g>\nonumber\\
&&~~~~~~+<\varphi-\varphi_{h},g>
+a_{h}(\varphi,\varphi^{*}-\varphi_{h}^{*})
-<f,\varphi^{*}-\varphi_{h}^{*}>+<f,\varphi^{*}-\varphi_{h}^{*}>\\
&&~~~=D_{h}^*(\varphi-\varphi_{h},\varphi^{*})+D_{h}(\varphi,\varphi^{*}-\varphi_{h}^{*})+<f,\varphi^{*}-\varphi_{h}^{*}>+<\varphi-\varphi_{h},g>.
\end{eqnarray*}
Combining the above two equalities we obtain
\begin{eqnarray*}
<\varphi-\varphi_{h},g>=a_{h}(\varphi-\varphi_{h},\varphi^{*}-\varphi_{h}^{*})
-D_{h}^*(\varphi-\varphi_{h},\varphi^{*})-D_{h}(\varphi,\varphi^{*}-\varphi_{h}^{*}),
\end{eqnarray*}
and we have by duality,
\begin{eqnarray*}
&&\|\varphi-\varphi_{h}\|_{0,\partial\Omega}=
\sup\limits_{g\in L^{2}(\partial\Omega)\setminus\{0\}}\frac{<\varphi-\varphi_{h},g>}{\|g\|_{0,\partial\Omega}}
\lesssim \|\varphi-\varphi_{h}\|_{h}\sup\limits_{g\in L^{2}(\partial\Omega)\setminus\{0\}}\frac{\|\varphi^{*}-\varphi_{h}^{*}\|_{h}}{\|g\|_{0,\partial\Omega}}\nonumber\\
&&~~~~~~+\sup\limits_{g\in L^{2}(\partial\Omega)\setminus\{0\}}\frac{|D_{h}(\varphi,\varphi^{*}-\varphi_{h}^{*})|+
|D_{h}^{*}(\varphi-\varphi_{h},\varphi^{*})|}{\|g\|_{0,\partial\Omega}}.
\end{eqnarray*}
Substituting (\ref{s3.7}), (\ref{s3.8}), (\ref{s3.25}) and (\ref{s3.27}) into the above inequality, and
using the regularity estimates (\ref{s2.24}) and (\ref{s2.28}),  we get (\ref{s3.26}).
Similarly we can prove (\ref{s3.28}). The proof is
completed.~~~$\Box$\\

\noindent{\bf Remark 2~} Consider the Neumann boundary problem:
find $\varphi\in H^{1}(\Omega)$ such that
\begin{eqnarray}\label{s3.29}
a(\varphi,v)= <f,v>+(\zeta,v),~~~\forall v\in H^{1}(\Omega).
\end{eqnarray}
 Let $\varphi$ and $\varphi_{h}$ be the exact solution and the C-R element solution of
(\ref{s3.29}), respectively,
\begin{eqnarray*}
D_{h}(\varphi,v)=a_{h}(\varphi,v)-<f, v>-(\zeta,v),
\end{eqnarray*}
and let $\varphi^{*}$ and $\varphi_{h}^{*}$ be the exact solution and the C-R element solution of
the dual problem of (\ref{s3.29}), respectively.
Then the analysis and conclusions in this section are also valid for (\ref{s3.29}).

\section{Error estimates of discrete Stekloff eigenvalues}
 \indent In this paper we suppose that
$\{\lambda_{j}\}$ and $\{\lambda_{j,h}\}$  are enumerations of the eigenvalues of (\ref{s2.1}) and (\ref{s2.10}) respectively according to the same sort rule,
and let $\lambda=\lambda_{m}$ be the $m$th eigenvalue with the
algebraic multiplicity $q$ and the ascent $\alpha$,
$\lambda_{m}=\lambda_{m+1}=\cdots,\lambda_{m+q-1}$. When
$\|T_{h}-T\|_{0,\partial\Omega}\to 0$, $q$ eigenvalues
$\lambda_{m,h},\cdots,\lambda_{m+q-1,h}$ of (\ref{s2.10}) will
converge to $\lambda$ (see Lemma 5 on page 1091 of \cite{dunford}).\\
 \indent Let $M(\lambda)$ be the space of generalized
eigenvectors associated with $\lambda$ and $T$, let $M_{h}(\lambda_{j,h})$
be the space of generalized eigenvectors associated with
$\lambda_{j,h}$ and $T_{h}$, and let
$M_{h}(\lambda)=\sum_{j=m}^{m+q-1}M_{h}(\lambda_{j,h})$.
 In view of the
dual problem (\ref{s2.6}) and (\ref{s2.14}), the definitions of
$M(\lambda^{*})$, $M_{h}(\lambda_{j,h}^{*})$ and
$M_{h}(\lambda^{*})$ are analogous to $M(\lambda)$, $M_{h}(\lambda_{j,h})$ and
$M_{h}(\lambda)$.\\
\indent Given two closed subspaces $V$ and $U$, denote
 \begin{eqnarray*}
\delta(V,U)=\sup\limits_{u\in V\atop\|u\|_{0,\partial\Omega}=1}\inf\limits_{v\in U}\|u-v\|_{0,\partial\Omega},
~~~\hat{\delta}(V,U)=\max\{\delta(V,U), \delta(U,V)\}.
\end{eqnarray*}
And denote
$
\hat{\lambda}_{h}=\frac{1}{q}\sum\limits_{j=m}^{m+q-1}\lambda_{j,h}
$. Thanks to \cite{babuska},  we get the following Theorem 4.\\
\noindent{\bf Theorem 4~}
 Suppose $M(\lambda), M(\lambda^{*})\subset H^{1+t}(\Omega)$~($t\in [r,1]$ for $\Omega\subset\mathbb{R}^2$, and $t\in [r3,1]$ for $\Omega\subset\mathbb{R}^3$), and {\bf Hypothesis $R(\Omega)$} holds. Then
\begin{eqnarray}\label{s4.1}
&&\hat{\delta}(M(\lambda),M_{h}(\lambda))\lesssim h^{s+t},\\\label{s4.2}
&&|\hat{\lambda}_{h}-\lambda|\lesssim h^{2t},\\\label{s4.3}
&&|\lambda-\lambda_{j,h}|\lesssim h^{\frac{2t}{\alpha}},~~~j=m,m+1,\cdots,m+q-1;
\end{eqnarray}
suppose $u_{h}$ is
an eigenfunction corresponding to $\lambda_{j,h}$
($j=m,m+1,\cdots,m+q-1$), $\|u_{h}\|_{0,\partial\Omega}=1$, then there exists an eigenfunction
$u$ corresponding to $\lambda$, such that
\begin{eqnarray}\label{s4.4}
&&\|u_{h}-u\|_{0,\partial\Omega}\lesssim  h^{(s+t)\frac{1}{\alpha}},\\\label{s4.5}
&&\|u_{h}-u\|_{h}\lesssim h^{t} + h^{(s+t)\frac{1}{\alpha}};
\end{eqnarray}
where s=$\frac{r}{2}$ when $\Omega\subset\mathbb{R}^2$, s=$r3$ when $\Omega\subset\mathbb{R}^3$.\\
\noindent{\bf Proof~}
Note that $\|Tf-T_{h}f\|_{0,\partial\Omega}=\|Af-A_{h}f\|_{0,\partial\Omega}=\|\varphi-\varphi_{h}\|_{0,\partial\Omega}$, from (\ref{s3.26}) with $t=s$ we deduce
\begin{eqnarray}\label{s4.6}
&&~~~\|T-T_{h}\|_{0,\partial\Omega}=\sup\limits_{f\in L^{2}(\partial\Omega),\|f\|_{0,\partial\Omega}=1 }\|Tf-T_{h}f\|_{0,\partial\Omega}\nonumber\\
&&\lesssim \sup\limits_{f\in L^{2}(\partial\Omega),\|f\|_{0,\partial\Omega}=1 } h^{2s}\|Af\|_{1+s}\lesssim  h^{2s}\|f\|_{0,\partial\Omega}\lesssim
h^{2s}\to 0~~~(h\to 0).
\end{eqnarray}
Thus from Theorem 7.1, Theorem 7.2  (inequality (7.12)), Theorem 7.3 and Theorem 7.4
in \cite{babuska} we get
\begin{eqnarray}\label{s4.7}
&&\hat{\delta}(M(\lambda),M_{h}(\lambda))) \lesssim \|(T-T_{h})\mid_{M(\lambda)}\|_{0,\partial\Omega},\\\label{s4.8}
&&\mid\lambda-\hat{\lambda}_{h}\mid \lesssim \sum\limits_{i,j=m}^{m+q-1}\mid<(T-T_{h})\varphi_{i},\varphi_{j}^*>\mid\nonumber\\
&&~~~+\|(T-T_{h}) \mid_{M(\lambda)} \|_{0,\partial\Omega}\|(T^*-T_{h}^*) \mid_{M(\lambda^{*})} \|_{0,\partial\Omega},\\\label{s4.9}
&&\mid\lambda-\lambda_{h}\mid \lesssim \{\sum\limits_{i,j=m}^{m+q-1}\mid<(T-T_{h})\varphi_{i},\varphi_{j}^*>\mid\nonumber\\
&&~~~+\|(T-T_{h}) \mid_{M(\lambda)} \|_{0,\partial\Omega}\|(T^*-T_{h}^*) \mid_{M(\lambda^{*})} \|_{0,\partial\Omega}\}^{1/\alpha},\\
\label{s4.10}
&&\|u_{h}-u\|_{0,\partial\Omega}\leq C\|(T_{h}-T)\mid_{M(\lambda)}\|_{0,\partial\Omega}^{\frac{1}{\alpha}},
\end{eqnarray}
where $\varphi_{m},\cdots,\varphi_{m+q-1}$ are any basis for $M(\lambda)$ and $\varphi_{m}^{*},\cdots,\varphi_{m+q-1}^{*}$ are the dual basis in $M(\lambda^{*})$.\\
From (\ref{s3.26}) with $t\in [r, 1]$ we obtain
\begin{eqnarray}\label{s4.11}
&&\|(T-T_{h})|_{M(\lambda)}\|_{0,\partial\Omega}=\sup\limits_{f\in M(\lambda),\|f\|_{0,\partial\Omega}=1 }\|Tf-T_{h}f\|_{0,\partial\Omega}\nonumber\\
&&~~~\lesssim h^{s+t} \sup\limits_{f\in M(\lambda),\|f\|_{0,\partial\Omega}=1 } \|Af\|_{1+t}.
\end{eqnarray}
Similarly we have
\begin{eqnarray}\label{s4.12}
&&\|(T^*-T_{h}^*)|_{M(\lambda^*)}\|_{0,\partial\Omega}\lesssim h^{s+t} \sup\limits_{f\in M(\lambda^*),\|f\|_{0,\partial\Omega}=1 } \|A^*f\|_{1+t}.
\end{eqnarray}
Substituting (\ref{s4.11}) into (\ref{s4.7}) and (\ref{s4.10}) we get (\ref{s4.1}) and (\ref{s4.4}), respectively.\\
\indent The remainder is to prove (\ref{s4.2}), (\ref{s4.3}) and (\ref{s4.5}). An easy calculation show that
\begin{eqnarray}\label{s4.13}
&&<(T-T_{h})\varphi_{i},\varphi_{j}^*>=<T\varphi_{i},\varphi_{j}^*>-<T_{h}\varphi_{i},\varphi_{j}^*>\nonumber\\
&&~~~=a_{h}(A\varphi_{i},A^*\varphi_{j}^*)-a_{h}(A_{h}\varphi_{i},A_{h}^*\varphi_{j}^*)\nonumber\\
&&~~~=a_{h}(A\varphi_{i}-A_{h}\varphi_{i},A^*\varphi_{j}^*)+a_{h}(A_{h}\varphi_{i},A^*\varphi_{j}^{*}-A_{h}^{*}\varphi_{j}^*)\nonumber\\
&&~~~=a_{h}(A\varphi_{i}-A_{h}\varphi_{i},A^*\varphi_{j}^*)+a_{h}(A\varphi_{i},A^*\varphi_{j}^*-A_{h}^*\varphi_{j}^*)\nonumber\\ &&~~~~~~-a_{h}(A\varphi_{i}-A_{h}\varphi_{i},A^*\varphi_{j}^*-A_{h}^*\varphi_{j}^*).
\end{eqnarray}
By (\ref{s3.1}) and (\ref{s3.2}) with $f=\varphi_{i}$,  $\varphi=A\varphi_{i}$, $g=\varphi_{j}^*$ and $\varphi^*=A^*\varphi_{j}^*$, we obtain
\begin{eqnarray*}
&&a_{h}(A\varphi_{i}-A_{h}\varphi_{i},A^*\varphi_{j}^*)=D_{h}^*(A\varphi_{i}-A_{h}\varphi_{i},A^*\varphi_{j}^*)+<T\varphi_{i}-T_{h}\varphi_{i},\varphi_{j}^*>,
\nonumber\\
&&a_{h}(A\varphi_{i},A^*\varphi_{j}^{*}-A_{h}^{*}\varphi_{j}^*)=D_{h}(A\varphi_{i},A^*\varphi_{j}^{*}-A_{h}^{*}\varphi_{j}^*)
+<\varphi_{i},T^*\varphi_{j}^{*}-T_{h}^*\varphi_{j}^*>
\nonumber\\
&&~~~=D_{h}(A\varphi_{i},A^*\varphi_{j}^{*}-A_{h}^{*}\varphi_{j}^*)
+<T\varphi_{i}-T_{h}\varphi_{i},\varphi_{j}^*>.
\end{eqnarray*}
Substituting the above two relations into (\ref{s4.13}) we get
\begin{eqnarray}\label{s4.14}
&&<(T-T_{h})\varphi_{i},\varphi_{j}^*>
=-D_{h}^{*}(A\varphi_{i}-A_{h}\varphi_{i},A^*\varphi_{j}^*)-D_{h}(A\varphi_{i},A^*\varphi_{j}^{*}-A_{h}^*\varphi_{j}^*)\nonumber\\ &&~~~~~~+a_{h}(A\varphi_{i}-A_{h}\varphi_{i},A^*\varphi_{j}-A_{h}\varphi_{j}^*),
\end{eqnarray}
which together with (\ref{s3.7}), (\ref{s3.8}), (\ref{s3.25}) and (\ref{s3.27}) yields
\begin{eqnarray}\label{s4.15}
|<(T-T_{h})\varphi_{i},\varphi_{j}^*>|\lesssim h^{2t}.
\end{eqnarray}
Substituting (\ref{s4.15}), (\ref{s4.11}) and (\ref{s4.12}) into (\ref{s4.8}) and (\ref{s4.9}) we get (\ref{s4.2}) and (\ref{s4.3}), respectively.\\
From (\ref{s2.1}) and (\ref{s2.4}) we get
\begin{eqnarray*}
a(u,v)=a(A(-\lambda u),v),~~~\forall v\in H^{1}(\Omega),
\end{eqnarray*}
noting that $k^2$ is not an eigenvalue
of (\ref{s2.3}), we have
$
u=-\lambda Au
$.
Similarly, using (\ref{s2.10}) and (\ref{s2.12}) we can get
$
u_{h}=-\lambda_{h} A_{h}u_{h}
$.
Thus from (\ref{s3.25}), (\ref{s2.18}), (\ref{s2.24}), (\ref{s2.28}), (\ref{s4.3}) and (\ref{s4.4}) we deduce
\begin{eqnarray*}
&&\|u_{h}+\lambda A_{h}u\|_{h}=\|-\lambda_{h}A_{h}u_{h}+\lambda A_{h}u\|_{h}=\|A_{h}(\lambda_{h}u_{h}-\lambda u)\|_{h}\nonumber\\
&&~~~\leq\|(A-A_{h})(\lambda u-\lambda_{h}u_{h})\|_{h}+\|A(\lambda_{h}u_{h}-\lambda u)\|_{h}\nonumber\\
&&~~~\lesssim h^{
s}\|A(\lambda u-\lambda_{h}u_{h})\|_{1+s}+\|\lambda_{h}u_{h}-\lambda u\|_{0,\partial\Omega}\lesssim\|\lambda_{h}u_{h}-\lambda u\|_{0,\partial\Omega}\lesssim  h^{(t+s)\frac{1}{\alpha}},
\end{eqnarray*}
and by the triangular inequality
\begin{eqnarray*}
\|u_{h}-u\|_{h}\leq\|u_{h}+\lambda A_{h}u\|_{h}+\|\lambda Au-\lambda A_{h}u\|_{h}\lesssim h^{(t+s)\frac{1}{\alpha}}+h^{t},
\end{eqnarray*}
i.e., (\ref{s4.5}) holds. The proof is completed.~~~$\Box$

\section{Numerical experiments}

Consider the problem (\ref{s1.1}) on the test domain $\Omega \subset \mathbb{R}^{2} $, where $\Omega=(-\frac{\sqrt{2}}{2}, \frac{\sqrt{2}}{2})^2$ is the square, or $\Omega=(-1, 1)^2\setminus([0, 1)\times (-1, 0])$ is an
L-shaped domain with the largest inner angle
$\omega=\frac{3}{2}\pi$, or $\Omega=(-\frac{\sqrt{2}}{2}, \frac{\sqrt{2}}{2})^2\setminus\{0\leq x\leq\frac{\sqrt{2}}{2}, y=0\}$
is the square with a slit which the largest inner angle $\omega=2\pi$, or $\Omega$ is the unit disk,
and $k=1$,
$n(x)=4$ or $n(x)=4+4i$.\\
\indent We use Matlab 2012a to solve (\ref{s1.1}) on a Lenovo ideaPad PC with 1.8GHZ CPU and 8GB RAM. Our program is compiled under the
package of iFEM \cite{chenl}.\\

Referring to \cite{cakoni1,liu}, when $n=4$ we sort eigenvalues in descending order,
and when $n=4+4i$ we arrange complex eigenvalues by their imaginary parts from large to small.

For the unit disk, the exact Stekloff eigenvalues are given by (5.3) in \cite{liu}, and when $n=4$ the largest six eigenvalues are
$$\lambda_{1}=5.151841,~~~\lambda_{2,3}=0.223578,~~\lambda_{4,5}=-1.269100,~~~\lambda_{6}=-2.472703,$$
and when $n=4+4i$ the four complex eigenvalues with the largest imaginary parts are
$$\lambda_{1}=-0.320506+3.121689i,~~~\lambda_{2,3}=-0.136861+1.396737i,~~\lambda_{4}=-1.353076+0.791723i.$$
For the $L$-shaped and the slit domain, the reference eigenvalues of the exact eigenvalues are listed in Tables 11-12.

\subsection{Numerical experiments on uniform meshes}

We adopt a uniform mesh $\pi_{h}$ for each domain.
The numerical results
 are listed in Tables 1-8.
The error curves of the C-R eigenvalues are showed in Figs. 1-4.

\begin{table}\small
  \caption{The eigenvalues on the square: $n=4$.}
    \begin{tabular}{ccccccc}
    \hline\noalign{\smallskip}
    ~dof~& $\lambda_{1,h}$ & $\lambda_{2,h}$ & $\lambda_{3,h}$ & $\lambda_{4,h}$&$\lambda_{5,h}$& $\lambda_{6,h}$ \\
    \noalign{\smallskip}\hline\noalign{\smallskip}
    3136  & 2.2018805  & -0.2116751  & -0.2116708  & -0.9069429  & -2.7589883  & -2.7522381  \\
    12416 & 2.2023533  & -0.2121076  & -0.2121070  & -0.9077740  & -2.7664177  & -2.7646187  \\
    49408 & 2.2024690  & -0.2122160  & -0.2122159  & -0.9079851  & -2.7683097  & -2.7678463  \\
    197120 & 2.2024977  & -0.2122431  & -0.2122431  & -0.9080383  & -2.7687870  & -2.7686695  \\
 \hline
    \end{tabular}%
\end{table}%
\begin{table}\small
  \caption{The eigenvalues on the L-shaped domain: $n=4$.}
    \begin{tabular}{ccccccc}
    \hline\noalign{\smallskip}
    ~dof~& $\lambda_{1,h}$ & $\lambda_{2,h}$ & $\lambda_{3,h}$ & $\lambda_{4,h}$&$\lambda_{5,h}$& $\lambda_{6,h}$ \\
    \noalign{\smallskip}\hline\noalign{\smallskip}
    9344  & 2.5335485  & 0.8592520  & 0.1246281  & -1.0845725  & -1.0901869  & -1.4147102  \\
    37120 & 2.5333019  & 0.8583814  & 0.1245509  & -1.0851154  & -1.0909141  & -1.4163502  \\
    147968 & 2.5332364  & 0.8580275  & 0.1245311  & -1.0852527  & -1.0911151  & -1.4167642  \\
    590848 & 2.5332194  & 0.8578847  & 0.1245261  & -1.0852873  & -1.0911726  & -1.4168682  \\
 \hline
    \end{tabular}%
\end{table}%
\begin{table}\small
  \caption{The eigenvalues on the square with a slit: $n=4$.}
    \begin{tabular}{ccccccc}
    \hline\noalign{\smallskip}
    ~dof~& $\lambda_{1,h}$ & $\lambda_{2,h}$ & $\lambda_{3,h}$ & $\lambda_{4,h}$&$\lambda_{5,h}$& $\lambda_{6,h}$ \\
    \noalign{\smallskip}\hline\noalign{\smallskip}
    12448 & 1.4848728  & 0.4698829  & -0.1840366  & -0.6898362  & -1.8987837  & -1.9264514  \\
    49472 & 1.4847611  & 0.4658257  & -0.1841411  & -0.6900139  & -1.8995947  & -1.9278655  \\
    197248 & 1.4847266  & 0.4637839  & -0.1841672  & -0.6900592  & -1.8998016  & -1.9283610  \\
    787712 & 1.4847163  & 0.4627589  & -0.1841737  & -0.6900708  & -1.8998539  & -1.9285538  \\
 \hline
    \end{tabular}%
\end{table}%

\begin{table}\small
  \caption{The eigenvalues on the unit disk: $n=4$.}
    \begin{tabular}{ccccccc}
    \hline\noalign{\smallskip}
    ~dof~& $\lambda_{1,h}$ & $\lambda_{2,h}$ & $\lambda_{3,h}$ & $\lambda_{4,h}$&$\lambda_{5,h}$& $\lambda_{6,h}$ \\
    \noalign{\smallskip}\hline\noalign{\smallskip}
128628 & 5.1514757  & 0.2235716  & 0.2235710  & -1.2689792  & -1.2689809  & -2.4724133  \\
201444 & 5.1516049  & 0.2235739  & 0.2235738  & -1.2690218  & -1.2690228  & -2.4725174  \\
359676 & 5.1517065  & 0.2235760  & 0.2235759  & -1.2690553  & -1.2690557  & -2.4725980  \\
809421 & 5.1517811  & 0.2235773  & 0.2235773  & -1.2690803  & -1.2690804  & -2.4726559  \\
\hline
    \end{tabular}%
\end{table}%

\begin{table}\small
  \caption{The eigenvalues on the square: $n=4+4i$.}
    \begin{tabular}{ccccccc}
    \hline\noalign{\smallskip}
    ~dof~& $\lambda_{1,h}$ & $\lambda_{2,h}$ & $\lambda_{3,h}$& $\lambda_{4,h}$ & $\lambda_{5,h}$ & $\lambda_{6,h}$ \\
    \noalign{\smallskip}\hline\noalign{\smallskip}
    3136  &0.687353 & -0.342514 &-0.342525 & -0.948908 &-2.779702& -2.786716  \\
          & +2.494448i & +0.85089i & +0.850899i & +0.539844i & +0.53745i & +0.539647i \\
    12416 & 0.686749 & -0.342915  & -0.342916 &-0.949807 & -2.792169 & -2.794033 \\
          & +2.495075i & +0.850782i & +0.850784i & +0.540029i & +0.539839i & +0.540444i \\
    49408 & 0.686601 &-0.343014 & -0.343014 &-0.950034 & -2.795417 & -2.795897 \\
          & +2.495238i & +0.850755i & +0.850756i & +0.540079i & +0.540498i & +0.540656i \\
    197120 & 0.686564 &-0.343038 & -0.343038 & -0.950091 & -2.796245 & -2.796367  \\
          & +2.495280i & +0.850749i & +0.850749i & +0.540092i & +0.540671i & +0.540711i \\
 \hline
    \end{tabular}%
\end{table}%

\begin{table}\small
  \caption{The eigenvalues on the L-shaped domain: $n=4+4i$.}
    \begin{tabular}{ccccccccc}
    \hline\noalign{\smallskip}
    ~dof~& $\lambda_{1,h}$ & $\lambda_{2,h}$ & $\lambda_{3,h}$& $\lambda_{4,h}$ & $\lambda_{5,h}$ & $\lambda_{6,h}$\\
    \noalign{\smallskip}\hline\noalign{\smallskip}
    9344  & 0.513857  & 0.398298  & -0.076964 &-1.438567  & -1.654555  & -2.513849 \\
          & +2.881404i & +1.459758i & +1.042587i & +0.803689i & +0.766423i & +0.570528i \\
    37120 &0.514176& 0.397512  & -0.077125  & -1.440022  &-1.656531 & -2.516699  \\
          & +2.882086i & +1.459328i & +1.042656i & +0.804437i & +0.766548i & +0.571289i \\
    147968 &0.514259 & 0.397218  & -0.077165  &-1.440388 & -1.657092  & -2.517426  \\
          & +2.882263i & +1.459129i & +1.042672i & +0.80463i & +0.766548i & +0.571486i \\
    590848 & 0.514280 & 0.397106  & -0.077175 & -1.440479  & -1.657258 & -2.517610 \\
          & +2.882308i & +1.459043i & +1.042677i & +0.804678i & +0.766534i & +0.571536i \\
 \hline
    \end{tabular}%
\end{table}%

\begin{table}\small
  \caption{The eigenvalues on the square with a slit: $n=4+4i$.}
    \begin{tabular}{ccccccc}
    \hline\noalign{\smallskip}
    ~dof~& $\lambda_{1,h}$ & $\lambda_{2,h}$ & $\lambda_{3,h}$ & $\lambda_{4,h}$&$\lambda_{5,h}$& $\lambda_{6,h}$ \\
    \noalign{\smallskip}\hline\noalign{\smallskip}
    12448 &0.918974  & 0.299813 & -0.262446  & -0.741837 & -2.615356 &-2.840331 \\
          & +1.770802i & +1.003519i & +0.757437i & +0.608741i & +0.561764i & +0.493956i \\
    49472 & 0.919206 & 0.296211  & -0.262573 &-0.742028  & -2.618344 & -2.845935 \\
          & +1.770795i & +1.001745i & +0.757447i & +0.608765i & +0.562409i & +0.493673i \\
    197248 & 0.919276 & 0.294417  & -0.262604 & -0.742076 & -2.619113 & -2.847993 \\
          & +1.770791i & +1.000826i & +0.757449i & +0.608772i & +0.562579i & +0.493444i \\
    787712 & 0.919297  & 0.293522  &-0.262612 &-0.742088 & -2.619310 & -2.848830 \\
          & +1.770789i & +1.000356i & +0.75745i & +0.608774i & +0.562623i & +0.493306i \\
 \hline
    \end{tabular}%
\end{table}%

\begin{table}\small
  \centering
  \caption{The eigenvalues on the unit disk: $n=4+4i$.}
    \begin{tabular}{ccccccc}
    \hline\noalign{\smallskip}
    ~dof~& $\lambda_{1,h}$ & $\lambda_{2,h}$ & $\lambda_{3,h}$& $\lambda_{4,h}$ \\
    \noalign{\smallskip}\hline\noalign{\smallskip}
128628 & -0.320420  & -0.136864  & -0.136865  & -1.352964  \\
      & +3.124755i &  +1.39673i &  +1.39673i &  +0.79174i \\
201444 & -0.320451  & -0.136863  & -0.136863  & -1.353004  \\
          & +3.124732i & +1.396733i & +1.396733i & +0.791734i \\
359676 & -0.320475  & -0.136862  & -0.136862  & -1.353035  \\
         & +3.124714i & +1.396735i & +1.396735i &  +0.79173i \\
 809421 & -0.320492  & -0.136862  & -0.136862  & -1.353058  \\
         &   +3.124700i & +1.396736i & +1.396736i & +0.791726i \\
 \hline
    \end{tabular}%
\end{table}%

\begin{figure}
  \centering
  \includegraphics[width=5.5cm,height=4.5cm]{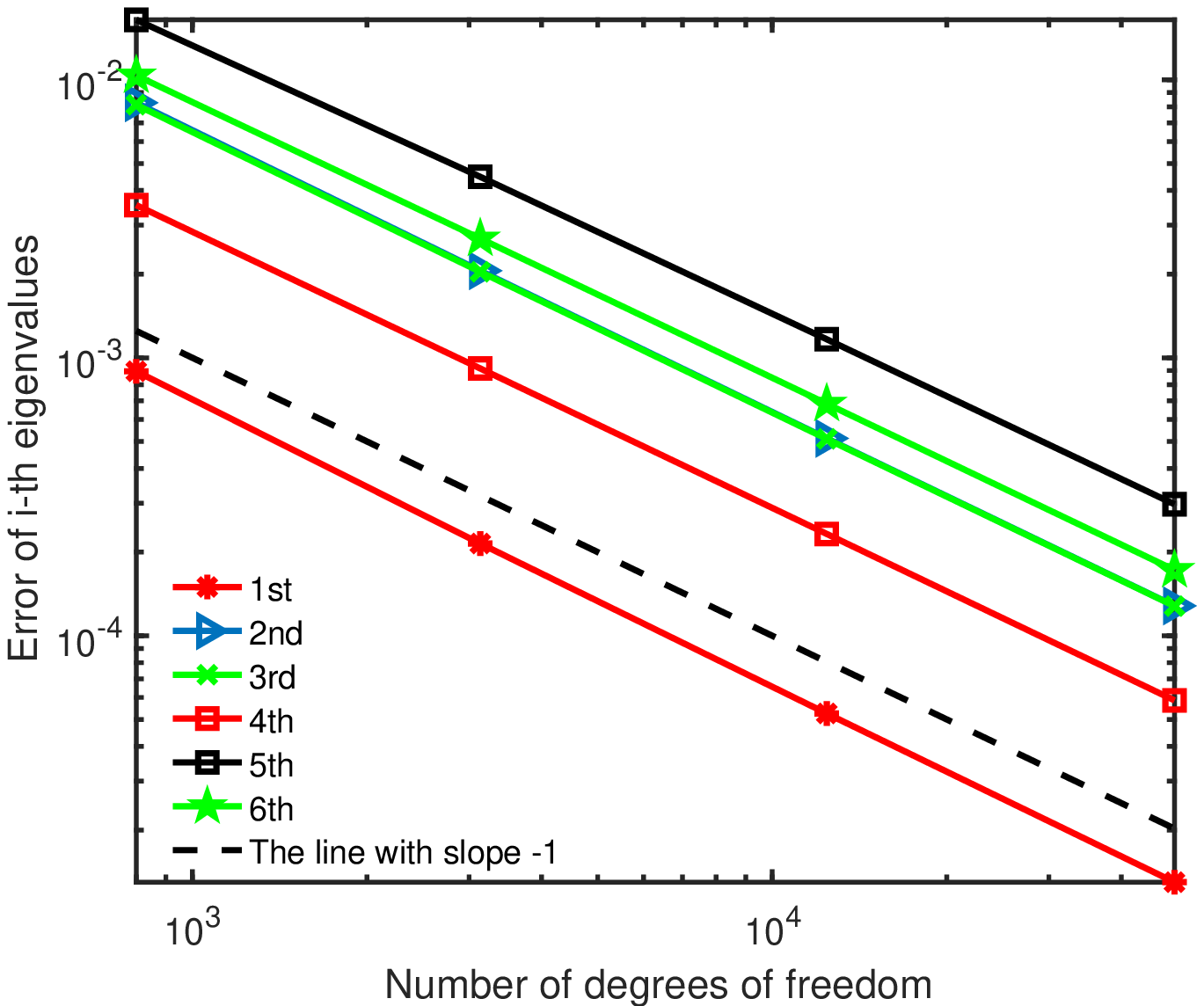}
  \includegraphics[width=5.5cm,height=4.5cm]{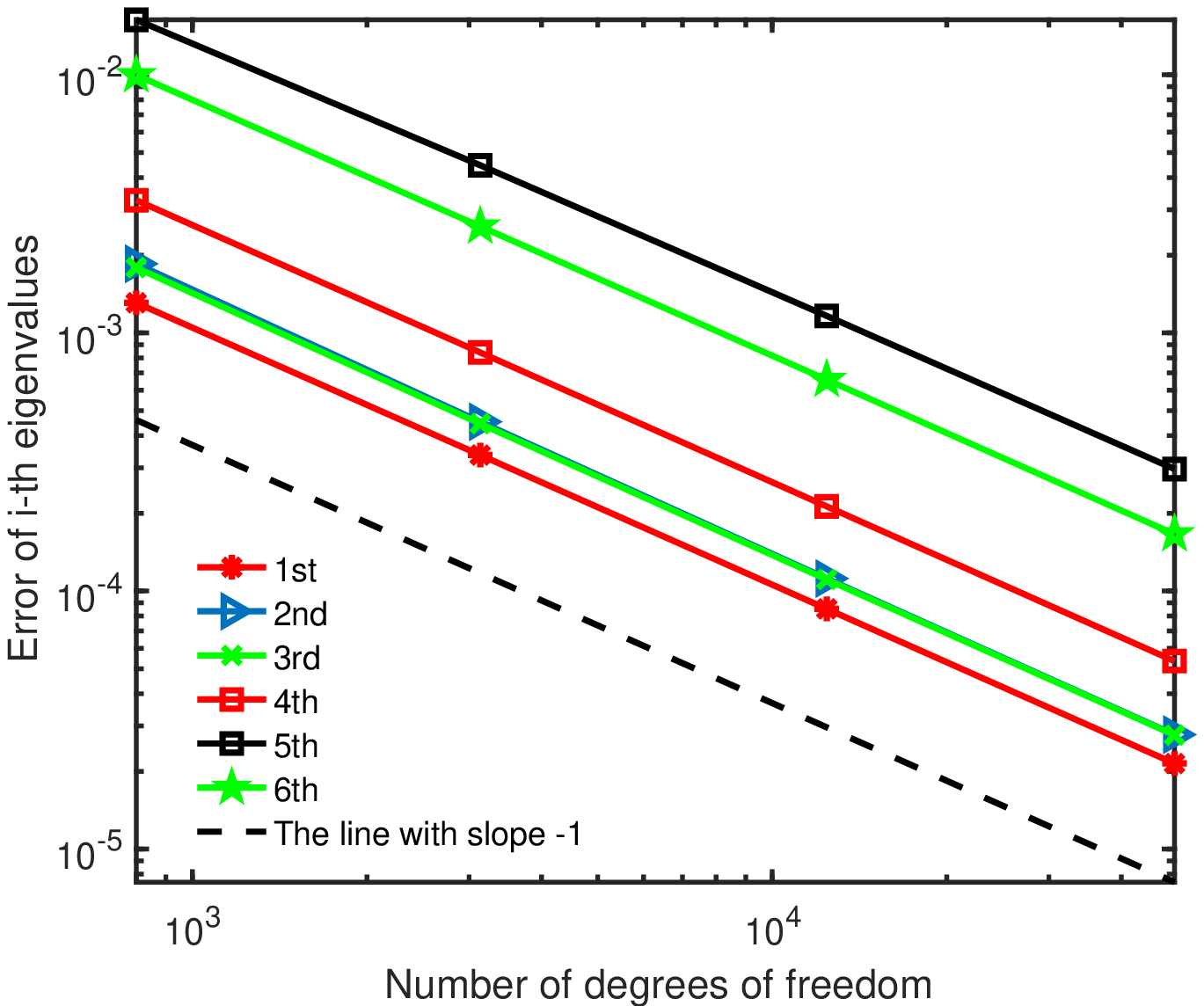}
  \caption{\emph{ The error curves of the first six eigenvalues on the square (left: $n=4$, right: $n=4+4i$)}}
  \centering
  \includegraphics[width=5.5cm,height=4.5cm]{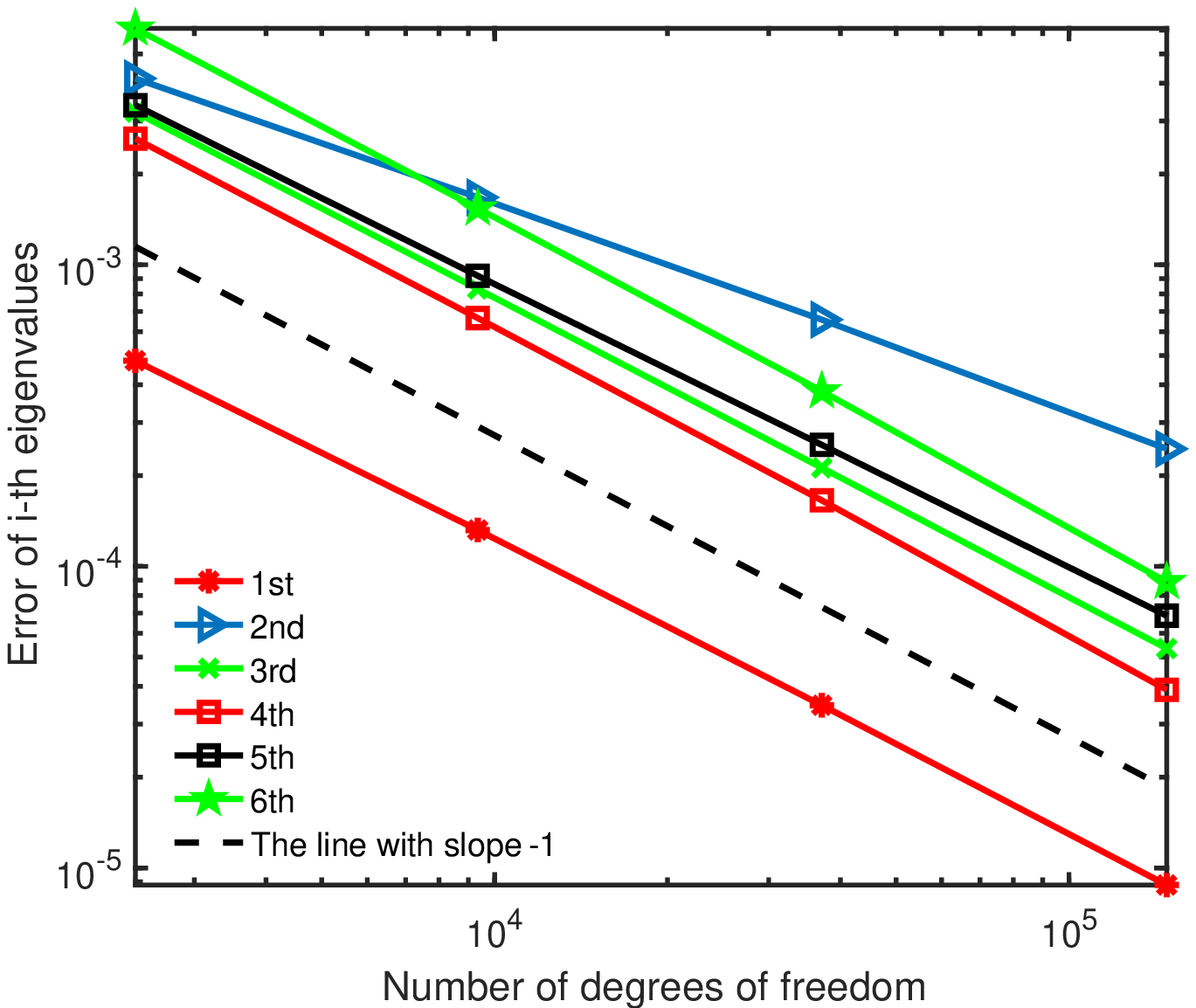}
  \includegraphics[width=5.5cm,height=4.5cm]{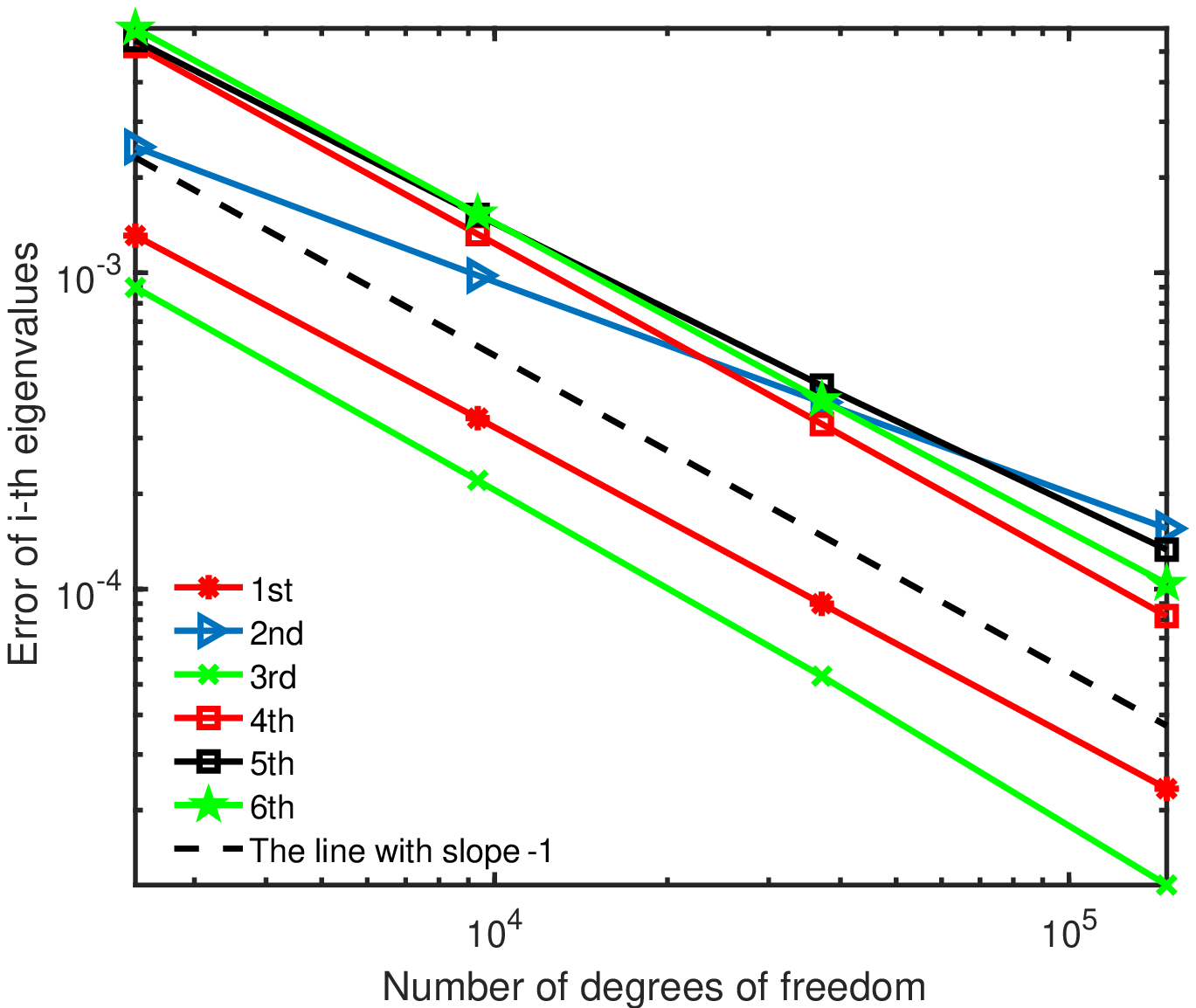}
  \caption{\emph{ The error curves of the first six eigenvalues on the L-shaped domain (left: $n=4$, right: $n=4+4i$)}}
  \centering
  \includegraphics[width=5.5cm,height=4.5cm]{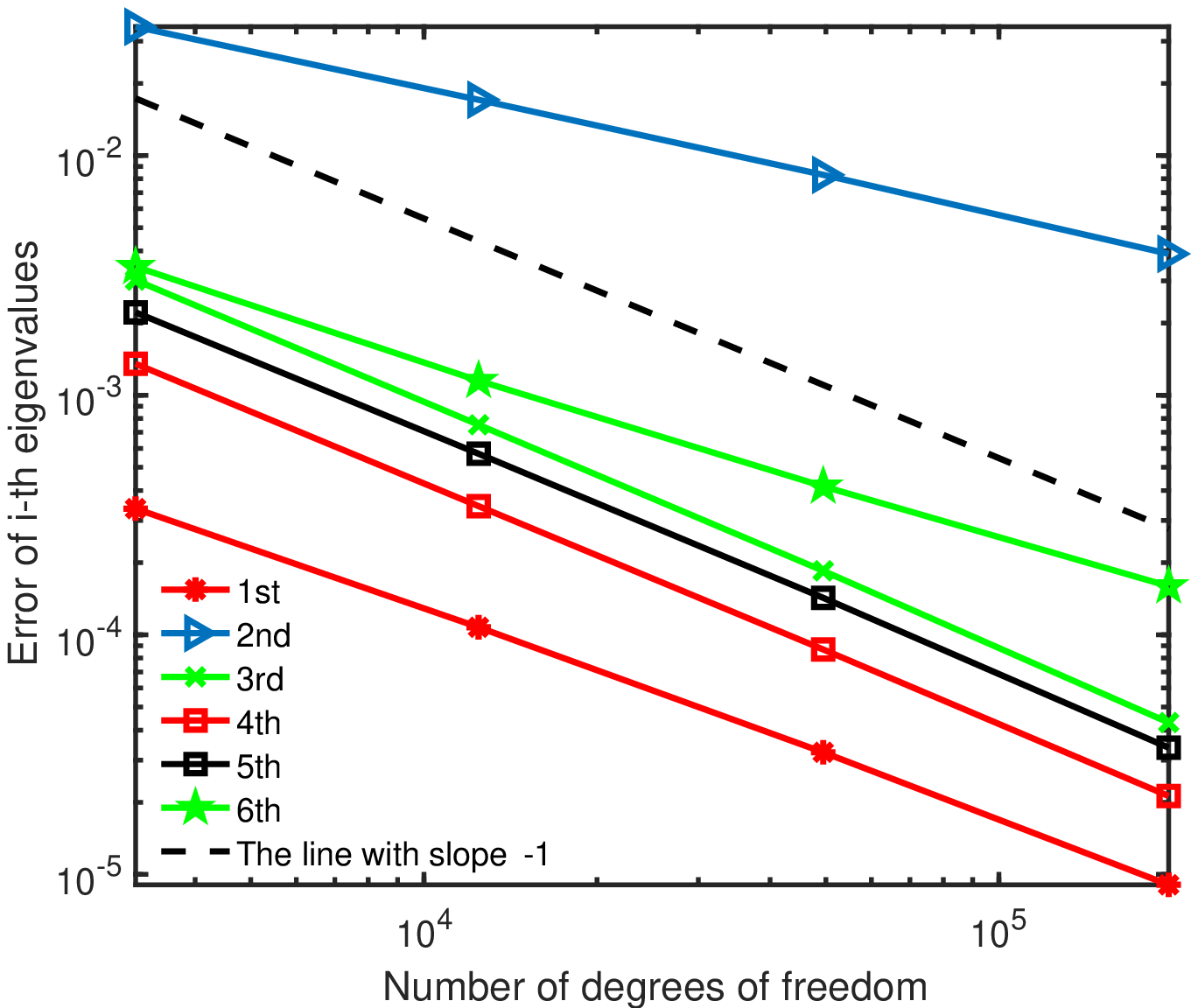}
  \includegraphics[width=5.5cm,height=4.5cm]{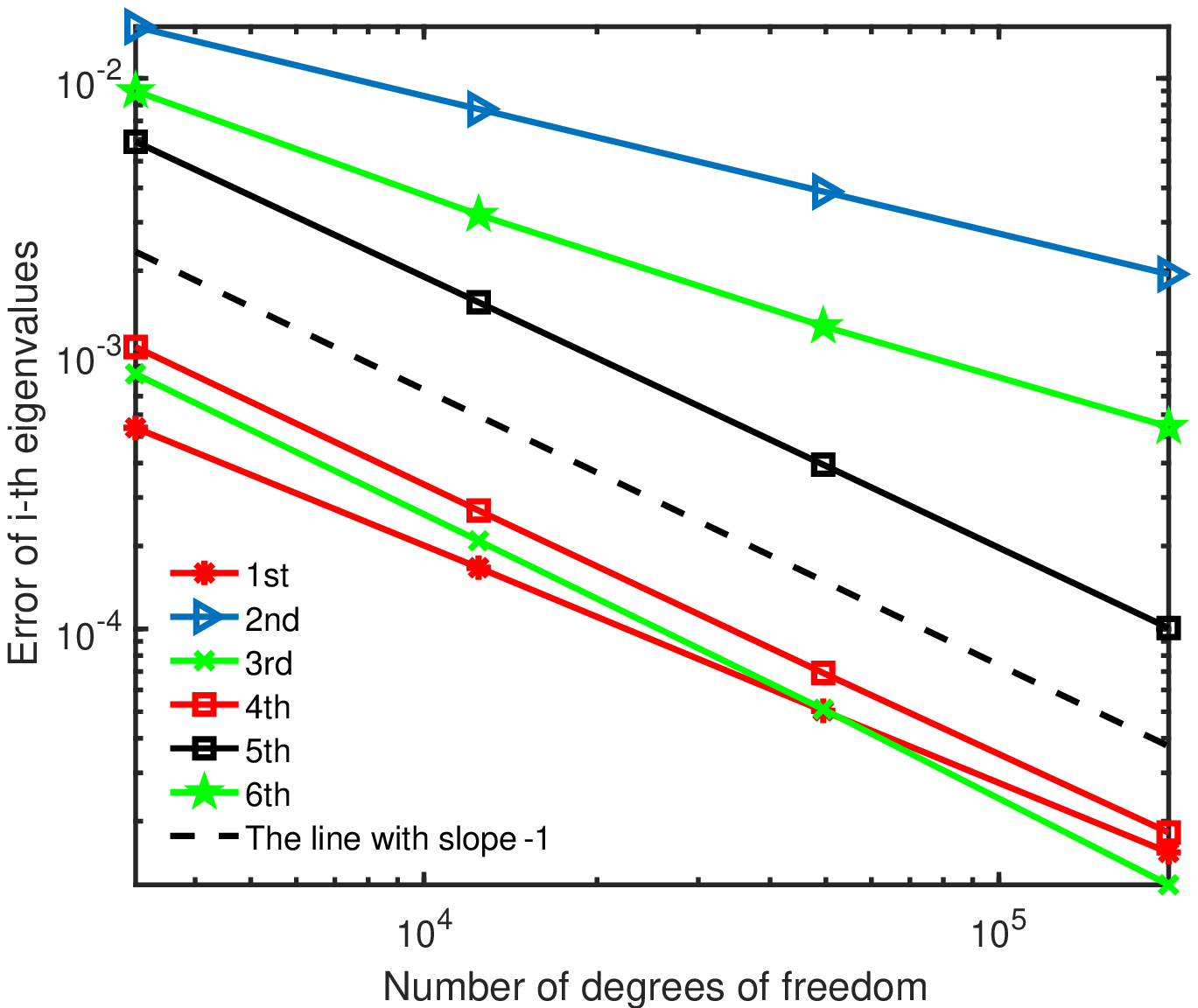}
  \caption{\emph{ The error curves of the first six eigenvalues on the square with a slit (left: $n=4$, right: $n=4+4i$)}}
\end{figure}

From Lemma 2, the regularity results, we know that for the square domain $2r=2$, for the L-shaped domain $2r\approx\frac{4}{3}$,
for the unit square with a slit $2r\approx1$.
From Fig. 1 we can see that the convergence order of $\lambda_{1,h},\lambda_{2,h}\cdots,\lambda_{6,h}$ are approximately equal to
$2$ on  the square domain; from Fig. 2 we can see that the convergence order of $\lambda_{2,h}$ is approximately equal to
$\frac{4}{3}\approx 1.333333$ on the L-shaped domain, and the eigenfunction corresponding to $\lambda_{2}$
has lower smoothness than others;
from Fig. 3 we can see that the convergence order of $\lambda_{2,h}$ is approximately equal to $1$ on the slit domain, and the eigenfunction corresponding to $\lambda_{2}$ is also less smoother that others, which are coincide with
the theoretical results.
Although there is an effect of reginal approximation for the computation on a disk, namely,
replacing the disk $\Omega$ with a similar polygonal $\Omega^{h}$, from Tables 7-8 and Fig. 4 we see that C-R element eigenvalues can approximate the exact ones.\\
For the square and the L-shaped domain, we also compare the numerical results in Tables 1-2, 5-6 with Tables 5.2-5.3, 5.5-5.6 in \cite{liu} and find that,
with the increase of dof (or the decrease of mesh size $h$), the eigenvalues obtained by C-R element and the conforming element are getting closer.

\subsection{Numerical experiments on adaptive meshes}
In practical finite element computations, it is desirable to carry out the
computations in an adaptive fashion (see, e.g.,\cite{ainsworth2,babuska1,brenner1,shiwang,verfurth}
 and references cited therein).
 For the C-R element approximation of
 Steklov eigenvalue problem, the a posteriori error estimates has been developed by
\cite{russo}.
Referring to \cite{russo} in this subsection we give the a posteriori error estimators by formal deduction, and implement adaptive computation for (\ref{s1.1}).\\
\indent Let $\ell\in\mathcal{E}_{h}^{i}$
shared by elements $\kappa_{1}$ and $\kappa_{2}$, i.e., $\ell=\partial\kappa_{1}\cap\partial\kappa_{2}$.
We choose a unit normal vector $\gamma_{\ell}$, pointing outwards $\kappa_{2}$, and we set
 the jumps of the normal derivatives of $v_{h}$ across $\ell$ as
\begin{eqnarray*}
[[\bigtriangledown v_{h}]]_{\gamma}=\bigtriangledown v_{h}|_{\kappa_{2}}\cdot\gamma_{\ell}-\bigtriangledown v_{h}|_{\kappa_{1}}\cdot \gamma_{\ell}.
\end{eqnarray*}
Denote $\gamma_{\ell}=(\gamma_{\ell1},\gamma_{\ell2})$, then
the tangent $t_{\ell}=(-\gamma_{\ell2},\gamma_{\ell1})$ on $\ell$, and we write
the jumps of the tangential derivatives of $v_{h}$ across $\ell$ as
\begin{eqnarray*}
[[\bigtriangledown v_{h}]]_{t}=\bigtriangledown v_{h}|_{\kappa_{2}}\cdot t_{\ell}-\bigtriangledown v_{h}|_{\kappa_{1}}\cdot t_{\ell}.
\end{eqnarray*}
 Notice that these values
are independent of the chosen direction of the normal vector $\gamma_{\ell}$.\\
\indent Now we define the a posteriori error indicators $\eta_{\kappa}(u_{h})$ on $\kappa$ and $\eta(u_{h})$
on $\Omega$ for the primal
eigenfunction $u_{h}$:\\
For each $\ell\in \mathcal{E}$, let
\begin{eqnarray*}
J_{\ell,t}(u_{h})= \left \{
\begin{array}{ll}
[[\nabla u_{h}]]_{t},&~if~\ell\in\mathcal{E}^{i},\\
0&~if~\ell\in\mathcal{E}^{b},
\end{array}
\right.~~~
J_{\ell,\gamma}(u_{h})= \left \{
\begin{array}{ll}
[[\nabla u_{h}]]_{\gamma},&~if~\ell\in\mathcal{E}^{i},\\
2(\nabla u_{h}\cdot\gamma_{\ell}-\lambda_{h}u_{h})|_{\ell}&~if~\ell\in\mathcal{E}^{b},
\end{array}
\right.
\end{eqnarray*}
and let
\begin{eqnarray*}
&&\eta_{\kappa}(u_{h})^2=|\kappa|\|k^{2}nu_{h}\|_{0,\kappa}^{2}+\frac{1}{2}\sum\limits_{\ell\in \partial\kappa}|\ell|\|J_{\ell,\gamma}(u_{h})\|_{0,\ell}^{2}+\frac{1}{2}\sum\limits_{\ell\in \partial\kappa}|\ell|\|J_{\ell,t}(u_{h})\|_{0,\ell}^{2},\\\label{s5.4}
&&\eta^2(u_{h})=\sum\limits_{\kappa\in\pi_{h}}\eta_{\kappa}(u_{h})^2.
\end{eqnarray*}
Similarly,  we define the  a posteriori error indicators $\eta_{\kappa}(u_{h}^*)$ on $\kappa$ and $\eta(u_{h}^*)$
on $\Omega$  for the dual eigenfunction $u_{h}^{*}$.\\
\indent We use $\sum\limits_{\kappa\in\pi_{h}}(\eta_\kappa^2(u_h)+\eta_\kappa^2(u_h^{*}))$ as
the a posteriori error indicator of
$\lambda_{h}$.
Using the indicator and consulting the existing standard adaptive algorithms (see, e.g., \cite{chenl,dai1,mao}), we solve (\ref{s1.1}).
From Figs. 2-3  we find that the eigenfunction associated with $\lambda_{2}$ is singular, so in our numerical
experiments we compute the approximation of the second eigenvalue $\lambda_{2}$,
and the numerical results on the L-shaped domain and the slit domain are listed in Table 9 and Table 10, respectively.\\
We show the curves of error and the a posteriori error estimators obtained by adaptive computing for the eigenvalue $\lambda_{2,h}$ in Figs. 5-6. It can be seen from them that the
error curves and the error estimators' curves are both basically parallel to the line with slope -1, which indicate that
the a posteriori error estimators of numerical eigenvalues are reliable and efficient and
$\lambda_{2,h}$ achieves the convergence rate $O(h^{2})$.
\begin{table}\small
  \caption{The second eigenvalues on adaptive meshes on the L-shaped domain.}
    \begin{tabular}{cccccc}
    \hline\noalign{\smallskip}
    $l$&dof & $\lambda_{2,h}  (n=4)$ &$l$&$dof$& $\lambda_{2,h} (n=4+4i)$ \\
    \noalign{\smallskip}\hline\noalign{\smallskip}
   1     & 9344   & 0.859246  & 1  & 9344  & 0.398302+1.459749i \\
   2     & 10022 & 0.858839  & 2  & 9494  & 0.398303+1.459755i \\
    25    & 202370 & 0.857844 &37&    216512 & 0.397113+1.459061i \\
    26    & 225490 & 0.857841 &38&    242302 & 0.397112+1.459059i \\
    27    & 249481 & 0.857838 &39&    286483 & 0.397108+1.459048i \\
    28    & 276807 & 0.857832 &40&    310425 & 0.397092+1.459023i \\
    29    & 331662 & 0.857821 &41&    340309 & 0.397092+1.459023i \\
    30    & 387329 & 0.857817 &42&    391833 & 0.397084+1.459016i \\
 \hline
    \end{tabular}%
\end{table}%
\begin{table}\small
  \caption{The second eigenvalues on adaptive meshes on the square with a slit.}
    \begin{tabular}{cccccc}
    \hline\noalign{\smallskip}
    $l$&dof & $\lambda_{2,h}  (n=4)$ &$l$& dof & $\lambda_{2,h} (n=4+4i)$ \\
  \noalign{\smallskip}\hline\noalign{\smallskip}
    1     & 12448 & 0.469884  & 1     & 12448 & 0.299812+1.003523i \\
    2     & 12472 & 0.467948  & 2     & 12607 & 0.298475+1.002875i \\
%
    25    & 170854 & 0.461819  &  60    & 241012 & 0.292782+0.999958i \\
    26    & 192640 & 0.461803  &  61    & 250930 & 0.292782+0.999958i \\
    27    & 222566 & 0.461788  &  62    & 260992 & 0.292761+0.999943i \\
    28    & 261309 & 0.461786  &  63    & 295455 & 0.292741+0.999931i \\
    29    & 298511 & 0.461783  &  64    & 311350 & 0.292738+0.999931i \\
    30    & 335598 & 0.461779  &  65    & 338930 & 0.292725+0.999926i \\

     \hline
    \end{tabular}%
\end{table}%

\begin{figure}
  \centering
  \includegraphics[width=5.5cm,height=4.5cm]{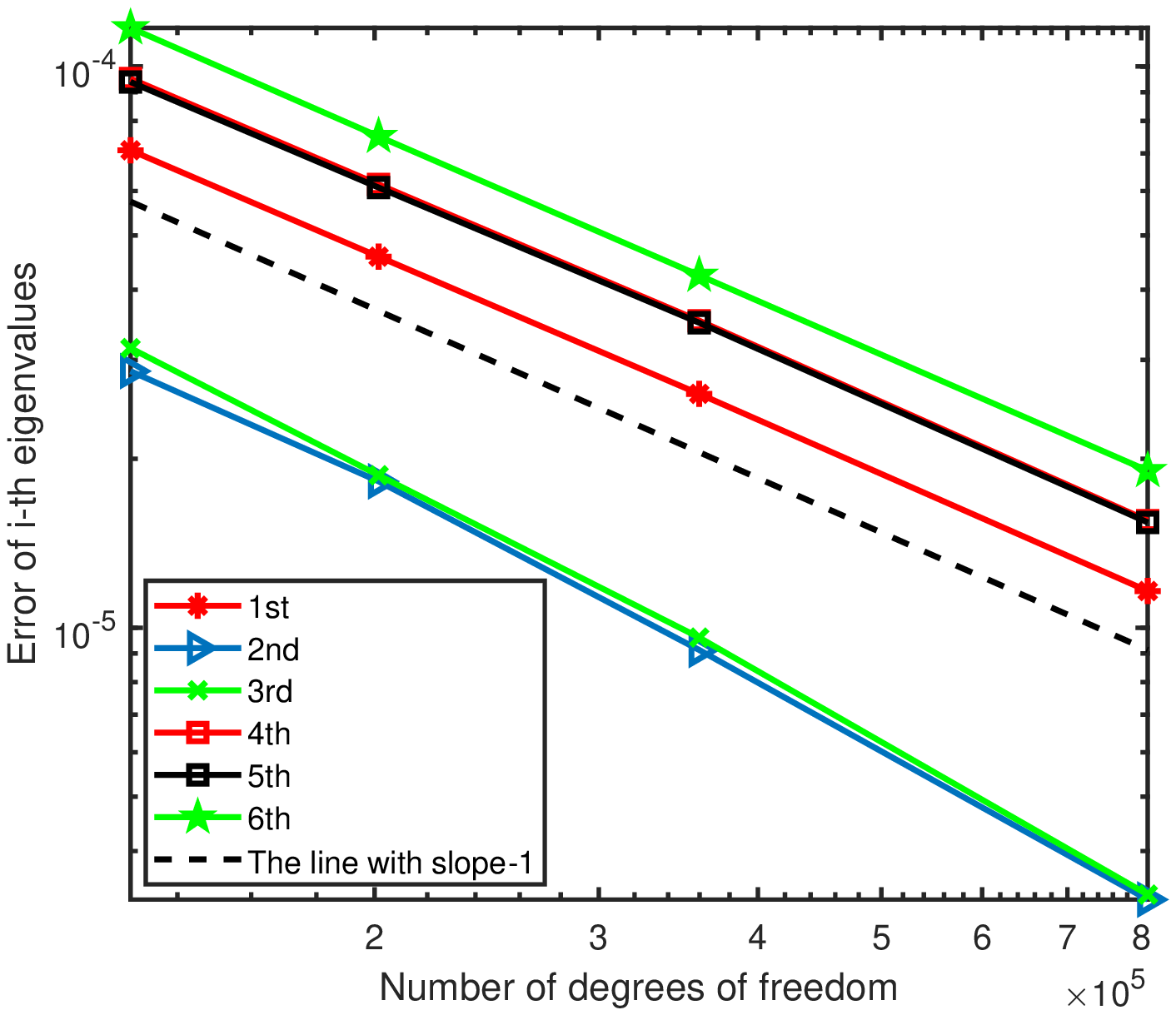}
  \includegraphics[width=5.5cm,height=4.5cm]{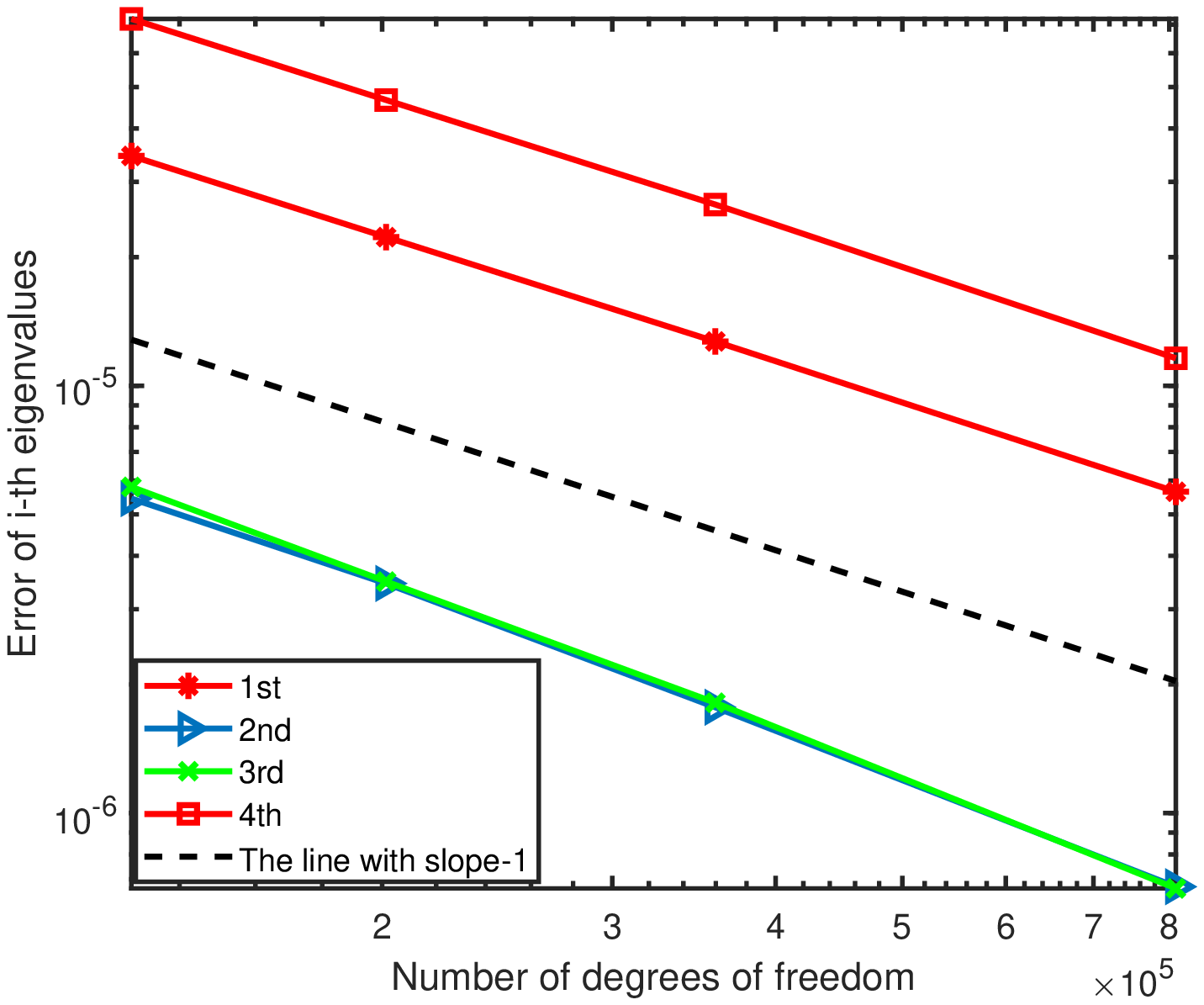}
  \caption{\emph{ The error curves of the eigenvalues for the unit disk on uniform meshes(left: $n=4$, right: $n=4+4i$)}}
  \centering
  \includegraphics[width=5.5cm,height=4.5cm]{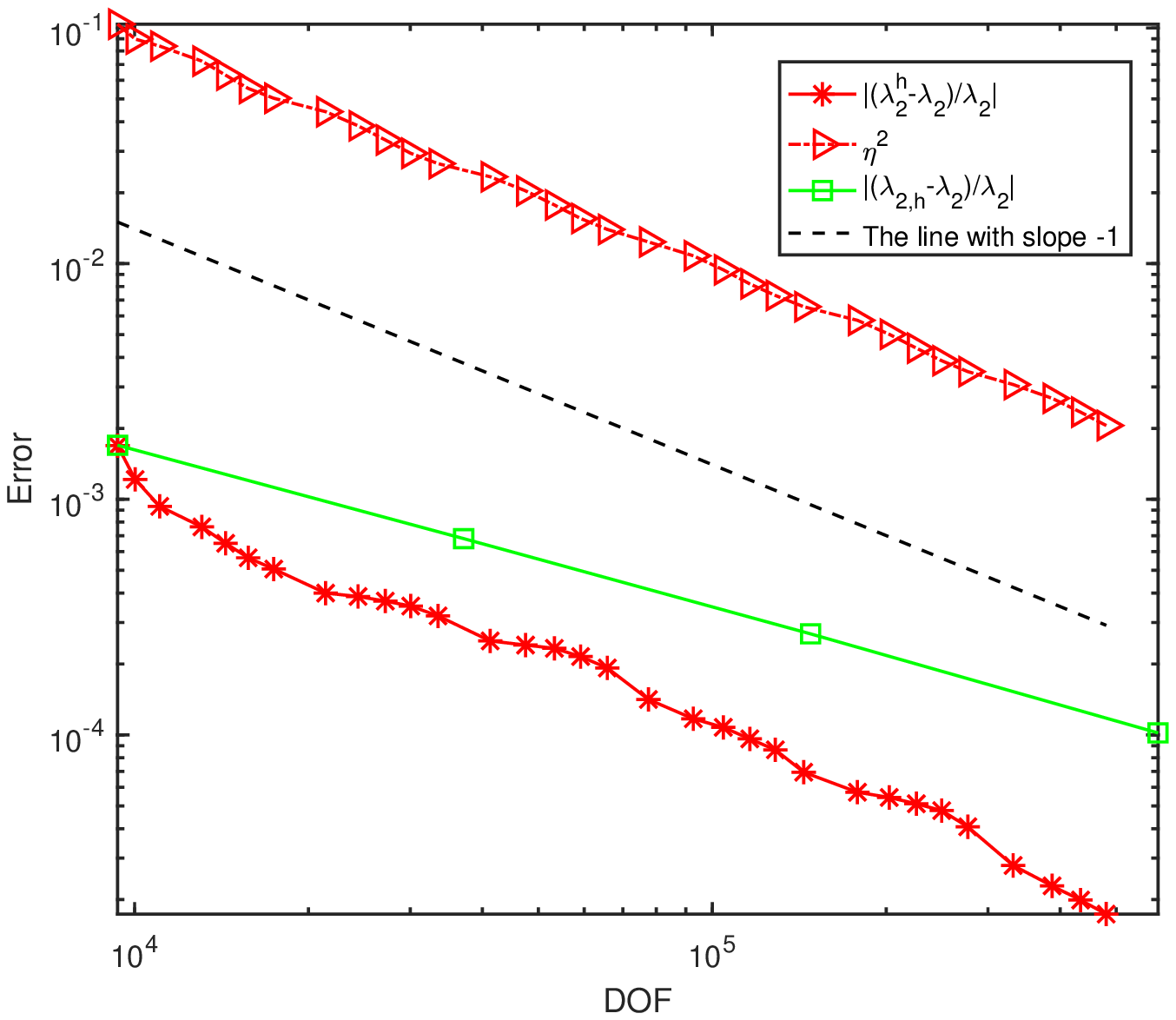}
  \includegraphics[width=5.5cm,height=4.5cm]{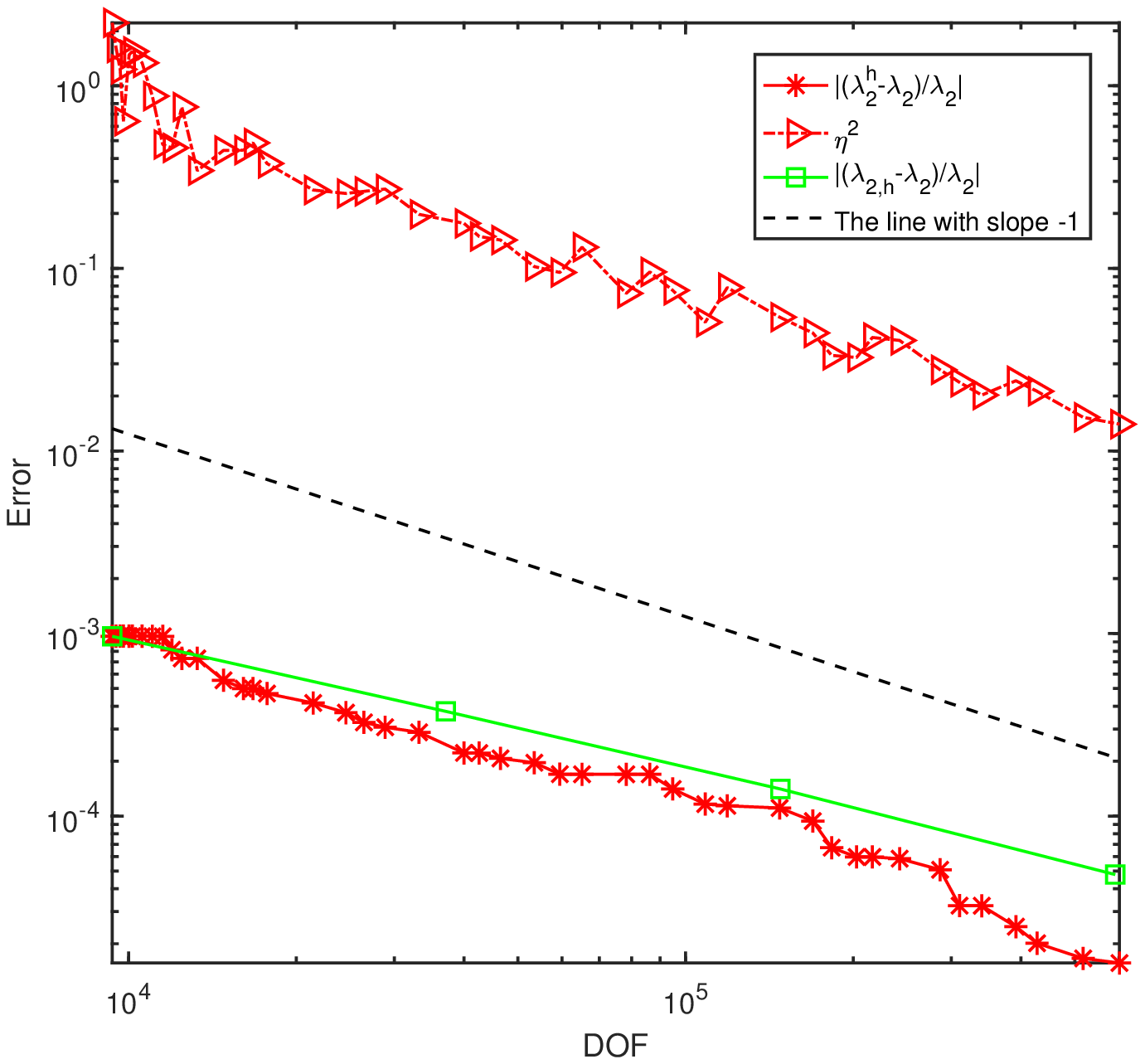}
  \caption{\emph{The error curves of the second eigenvalues on the L-shaped domain (left: $n=4$, right: $n=4+4i$)}}
  \centering
  \includegraphics[width=5.5cm,height=4.5cm]{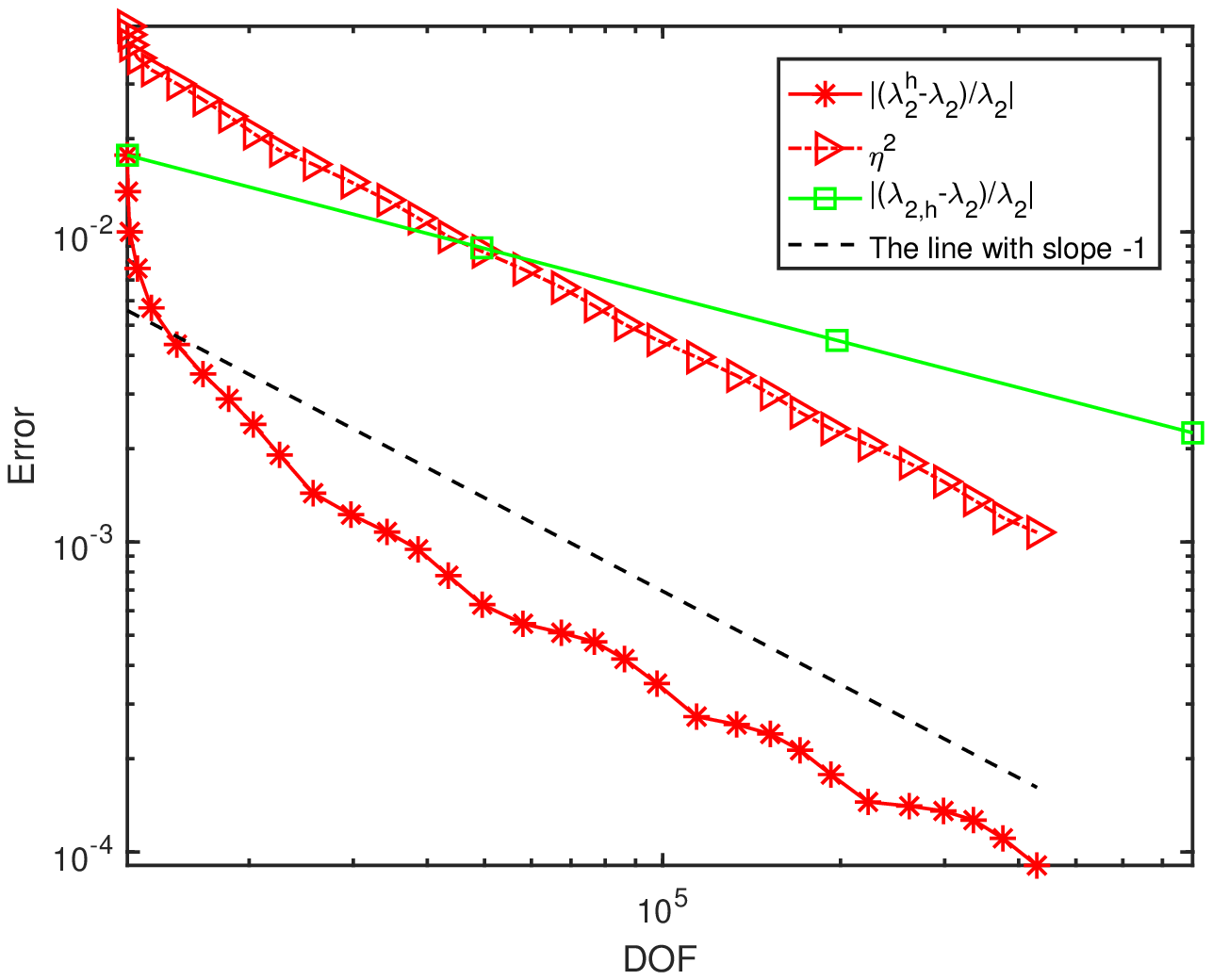}
  \includegraphics[width=5.5cm,height=4.5cm]{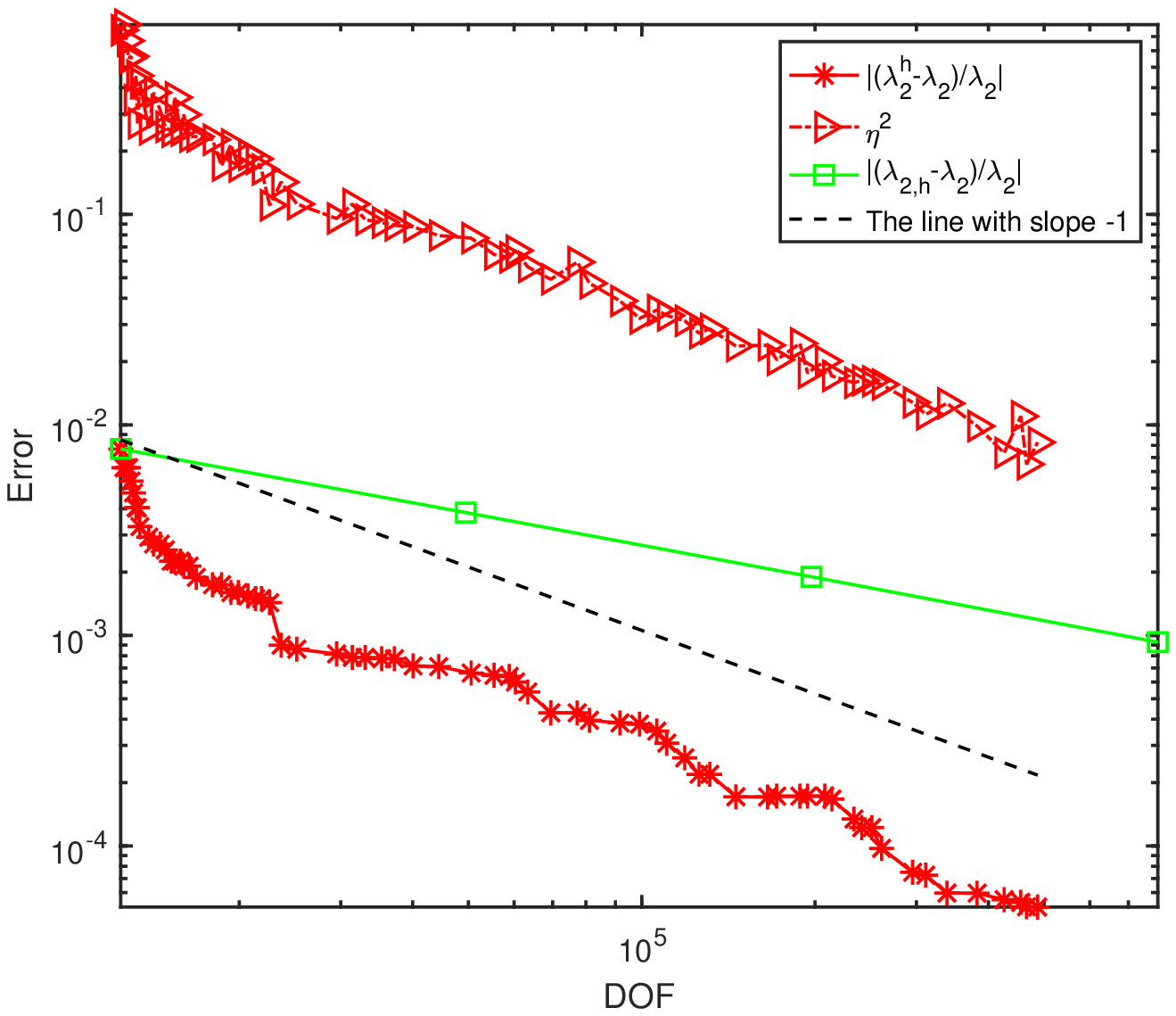}
  \caption{\emph{ The error curves of the second eigenvalues on the square with a slit (left: $n=4$, right: $n=4+4i$)}}
\end{figure}
From tables and figures we also see that under the same dof, the accuracy of approximate eigenvalues computed on adaptive meshes is far higher than
that of approximate eigenvalues computed on uniform meshes.\\

\noindent{\bf Remark 3}(The lower/upper bound of the Stekloff eigenvalues).\\
\indent We find in Tables 1-4 and 9-10 that when the index of refraction $n(x)$ is real, all series
of eigenvalues computed by C-R element show the tendency to decrease as the increase of $dof$
except the first one $\lambda_{1}$ on
the square and $\lambda_{1}\sim\lambda_{3}$ on the unit disk.
Note that the numerical results in \cite{cakoni1,liu} indicate that the conforming finite element eigenvalues approximate the exact ones from below when $n(x)$ is real.
So we also use the P1 conforming element to compute, and
obtain reference values of the exact eigenvalues by averaging the P1 conforming eigenvalues $\lambda_{j,h}^C$ and the C-R element eigenvalues $\lambda_{j,h}$. We list them in Tables 11-12.
The property of monotone convergence of the conforming finite element eigenvalues is easy to prove.
However, the property of monotone convergence of the C-R noconforming finite element eigenvalues is a non-trivial result. In 2014, Carstensen
and Gedicke \cite{carstensen1} prove rigorously the monotonicity for the classical Laplace eigenvalue problem. For the Stekloff eigenvalue problem considered in this paper, it is meaningful to study the monotone convergence of C-R element eigenvalues which is our next goal.\\
\begin{table}\small
  \caption{The reference eigenvalues $\lambda_j$(L) on the L-shaped domain and $\lambda_j$(Slit) on the square with a slit: $n=4$.}
    \begin{tabular}{ccccccc}
    \hline\noalign{\smallskip}
   $j$&$\lambda_{j,h}$& $\lambda_{j,h}^{C}$& $\lambda_j$(L) &$\lambda_{j,h}$&$\lambda_{j,h}^{C}$ &$\lambda_j$(Slit) \\
    \noalign{\smallskip}\hline\noalign{\smallskip}
    1&2.533219&2.533209&2.533214&1.484716&1.484710&1.484713\\	
    2&0.8578847&0.8577495&0.8578171&0.4627589&0.4612150&0.4619870\\
    3&0.1245261&0.1245229&0.1245245&-0.1841737&-0.1841765&-0.1841751\\
    4&-1.085287&-1.085303&-1.085295&-0.6900708&-0.6900769&-0.6900738\\
    5&-1.091173&-1.091207&-1.091190&-1.899854&-1.899878&-1.899866\\
    6&-1.416868&-1.416912&-1.416890&-1.928554&-1.928784&-1.928669\\
    \hline
    \end{tabular}%
\end{table}%

\begin{table}\small
  \caption{The reference eigenvalues $\lambda_j$(L) on the L-shaped domain and $\lambda_j$(Slit) on the square with a slit:~~$n=4+4i$.}
    \begin{tabular}{ccccccc}
    \hline\noalign{\smallskip}
   $j$&$\lambda_{j,h}$ & $\lambda_{j,h}^{C}$ & $\lambda_j$(L)&$\lambda_{j,h}$ & $\lambda_{j,h}^{C}$ & $\lambda_j$(Slit)  \\
    \noalign{\smallskip}\hline\noalign{\smallskip}
   1&0.5142799 &0.5143106 &0.5142952 &0.9192965 & 0.9193164&0.9193065\\
   &+2.882308i&+2.882334i&+2.882321i&+1.770789i&+1.770782i&+1.770786i\\
   2&0.3971057&0.3969716&0.3970387&0.2935223&0.2917372&0.2926298\\
   &+1.459043i&+1.458911i&+1.458977i&+1.000356i&+0.9993946i&+0.9998754i\\
   	3&-0.0771754&-0.0771792&-0.0771773&-0.2626120&-0.2626151&-0.2626135  \\
   &+1.042677i&+1.042673i&+1.042675i&+0.7574501i&+0.7574481i&+0.7574491i\\
   	4&-1.440479&-1.440535&-1.440507&-0.7420884&-0.7420981&-0.7420933\\
   &+0.8046784i&+0.8047093i&+0.8046939i&+0.6087744i&+0.6087755i&+0.6087749i\\
   	5&-1.657258&-1.657409&-1.657333&-2.619310 & -2.619442&-2.619376\\
   &+0.7665341i&+0.7664933i&+0.7665137i&+0.5626232i&+0.5626581i&+0.5626407i\\
   	6&-2.517610&-2.517765&-2.517687&-2.848830&-2.850239&-2.849534\\
   &+0.5715365i&+0.5715597i&+0.5715481i&+0.4933060i&+0.4929997i&+0.4931528i\\
    \hline
    \end{tabular}%
\end{table}%

\noindent {\bf Acknowledgments~}This work is supported by the National Natural Science Foundation of China
 (Grant Nos.11561014,11761022 ).


\begin{thebibliography}{s10}

\bibitem{ainsworth2}
Ainsworth, M., Oden, J.T.: A posteriori error estimates in the
finite element analysis. Wiley-Inter science, New York (2011)

\bibitem{ainsworth}Ainsworth, M.:  Robust a posteriori error estimation for nonconforming finite element approximation. SIAM J. Numer. Anal. 42, 2320-2341
(2005)

\bibitem{alonso}Alonso, A., Russo, A.D.:  Spectral approximation of
variationally-posed eigenvalue problems by nonconforming methods.
J. Comput. Appl. Math. 223, 177-197 (2009)

\bibitem{andreev}Andreev, A.B., Todorov, T.D.:  Isoparametric finite element approximation of a Steklov eigenvalue
problem. IMA. J. Numer. Anal. 24, 309-322 (2004)

\bibitem{armentano1}Armentano, M.G.: The effect of reduced integration in the Steklov eigenvalue problem.
Math. Mod. and Numer. Anal. $(M^{2}AN)$ 38, 27-36 (2004)

\bibitem{armentano2}Armentano, M.G., Padra, C.: A posteriori error estimates for the Steklov eigenvalue
problem. Appl. Numer. Math. 58, 593-601 (2008)

\bibitem{armentano3}Armentano, M.G., Duran, R.G.: Asymptotic lower bounds for eigenvalues by nonconforming finit element
methods. Electron. Trans. Numer. Anal. 17, 92-101 (2004)

\bibitem{babuska}Babuska, I., Osborn, J.E.: Eigenvalue Problems. in: P. G. Ciarlet, J. L. Lions(Eds),
Finite Element Methods (Part I), pp. 641-787, in: Handbook of Numerical Analysis,
Vol. 2, Elsevier Science Publishers, North-Holand (1991)

\bibitem{babuska1} Babuska, I., Rheinboldt, W.C.:  Error estimates for adaptive finite element
computations. SIAM J. Numer. Anal. 15, 736-754 (1978)

\bibitem{bergman}Bergman, S., Schiffer, M.:  Kernel Functions and Elliptic
Differential Equations in Mathematical Physics. Academic
Press, New York (1953)

\bibitem{bermudez} Bermudez, A.,  Rodriguez, R., Santamarina, D.: A finite element solution of an added mass formulation for coupled fluid-solid vibrations. Numer. Math. 87, 201-227 (2000)

\bibitem{bernardi}Bernardi, C., Hecht, F.: Error indicators for the mortar finite element discretization of Laplace equation.
Math. Comp. 71(240), 1371-1403 (2001)

\bibitem{bi} Bi, H., Yang, Y.: A two-grid method of the non-conforming Crouzeix-Raviart element for the Steklov
eigenvalue problem. Appl. Math. Comput. 217, 9669-9678 (2011)

\bibitem{boffi}Boffi, D.: Finite element approximation of eigenvalue
problems. Acta Numerica, 1-120 (2010)

 \bibitem{bramble}Bramble, J.H., Osborn, J.E.:  Approximation of Steklov eigenvalues of non-selfadjoint second order elliptic operators. in: A. K. Aziz, (Ed.), Math.Foundations of the Finite Element Method with Applications to PDE, PP.387-408, Academic, New York (1972)

\bibitem{brenner2}Brenner, S.C., Sung, L.Y.: Linear finite element methods for planar linear elasticity. Math. Comp. 59,
321-338 (1992)

\bibitem{brenner}Brenner, S.C., Scott, L.R.: The Mathematical Theory of Finite Element Methods. 2nd
ed.. Springer-Verlag, New york (2002)

\bibitem{brenner1}Brenner, S.C.: $C^{0}$ interior penalty methods. In Frontiers in Numerical Analysis-Durham 2010,
Lecture Notes in Computational Science and Engineering 85, pp.
79-147, Springer-Verlag (2012)

\bibitem{bucur}Bucur, D., Ionescu, I.R.: Asymptotic analysis and scaling of friction
parameters. Z. Angew. Math. Phys. (ZAMP) 57, 1042-1056 (2006)

\bibitem{caiz}Cai, Z., Ye, X., Zhang, S.: Discontinuous Galerkin finite element methods for interface problems: a priori and a posteriori
error estimations.
SIAM J. Numer. Anal. 49, 1761-1787 (2011)

\bibitem{cakoni1}Cakoni, F., Colton, D., Meng, S., Monk, P.: Stekloff eigenvalues in inverse scattering. SIAM
J. Appl. Math. 76(4), 1737-1763 (2016)

\bibitem{cao}Cao, L., Zhang, L., Allegretto, W., Lin, Y.: Multiscale asymptotic method for Steklov eigenvalue
equations in composite media. SIAM J. Numer. Anal. 51, 273-296 (2013)

\bibitem{carstensen5} Carstensen, C., Hu, J., Orlando, A.: Framework for the a posteriori error analysis of nonconforming finite elements. SIAM J. Numer.
Anal. 45, 68-82 (2007)

\bibitem{carstensen6}Carstensen, C., Hoppe, R.H.W.: Convergence analysis of an adaptive nonconforming finite element method. Numer. Math. 103,
 251-266 (2006)

\bibitem{carstensen1}Carstensen, C., Gedicke, J.: Guaranteed lower bounds for  eigenvalues. Math. Comp. 83, 2605-2629 (2014)

\bibitem{chatelin}Chatelin, F.: Spectral Approximations of Linear
Operators. Academic Press, New York, (1983)

\bibitem{chenl}
 Chen, L.: iFEM: an innovative finite element methods package in
MATLAB. Technical Report, University of California at Irvine (2009)

\bibitem{ciarlet}Ciarlet, P.G.:  Basic error estimates for elliptic
proplems. in: P. G. Ciarlet, J. L. Lions(Eds), Finite Element
Methods (Part I),  pp.21-343, in: Handbook of Numerical Analysis, Vol. 2,
Elsevier Science Publishers, North-Holand (1991)

\bibitem{conca}Conca, C., Planchard, J., Vanninathanm, M.: Fluid and Periodic
Structures. John Wiley \& Sons, New York (1995)

 \bibitem{crouzeix}Crouzeix, M., Raviart, P.A.: Conforming and nonconforming finite element methods for
solving the stationary stokes equations. RAIRO. Anal. Numer.
3, 33-75 (1973)

\bibitem{dai1}Dai, X., Xu, J., Zhou, A.: Convergence and optimal complexity of
adaptive finite element eigenvalue computations. Numer. Math. 110, 313-355 (2008)

\bibitem{dauge}Dauge, M.: Elliptic boundary value problems on corner domains:
smoothness and asymptotics of solutions. in: Lecture Notes in
Mathematics, vol.1341. Springer, Berlin (1988)

\bibitem{duran}Dari, E., Dur$\acute{a}$n, R., Padra, C., Vampa, V.: A posteriori error estimators for noconforming finite element methods.
 RAIRO Model. Math. Anal.
Numer. 30, 385-400 (1996)

\bibitem{dunford}Dunford, N., Schwartz, J. T.: Linear Operators, Vol.2: Spectral Theory,
Selfadjoint Operators in Hilbert Space. Interscience, New York,
(1963)

\bibitem{falk}Falk, R.S.: Nonconforming finite element methods for the equations of linear elasticity. Math. Comp.
57, 529-550 (1991)

\bibitem{garau}Garau, E.M., Morin, P.: Convergence and quasi-optimality of adaptive FEM for Steklov eigenvalue problems. IMA J. Numer. Anal.
 31(3), 914-946 (2011)

 \bibitem{girault}Girault, V., Raviart, P.A.: Finite Element Approximation of the Navier-Stokes
  Equations. Lecture Notes in Mathematics 749, Springer-Verlag, Berlin Heidelberg, New York, (1981)

\bibitem{huj}Hu, J., Huang, Y., Lin, Q.: The lower bounds for eigenvalues of elliptic operators by Nonconforming
finite element methods. J. Sci Comput. 61, 196-221 (2014)

\bibitem{jerison}Jerison, D.S., Kenig, C.E.: The Neumann problem on Lipschitz domains, Bull. Amer. Math. Soc. 4, 203-207 (1981)

\bibitem{kufner}Kufner, A., John, O., Fu$\check{c}$ik, S.: Function Spaces, Academia Publishing House, Prague (1977)

\bibitem{liq2}Li, Q., Lin, Q., Xie, H.: Nonconforming finite element approximations of the Steklov eigenvalue problems
and its lower bound approximations. Appl. Math. 58, 129-151 (2013)

\bibitem{lim}Li, M., Lin, Q., Zhang, S.:  Extrapolation and superconvergence of the Steklov eigenvalue
problems. Adv. Comput. Math. 33, 25-44 (2010)

\bibitem{liu}Liu, J., Sun, J., Turner, T.:
Spectral indicator method for a non-selfadjoint Steklov eigenvalue problem.
 ar Xiv: 1804.02582V1 [math. NA] 7 Apr (2018)


\bibitem{mao}Morin, P., Nochetto, R.H., Siebert, K.: Convergence of adaptive finite element methods.
SIAM Rev. 44, 631-658 (2002)

\bibitem{oden}Oden, J.T., Reddy, J.N.:  An Introduction to the Mathematical Theory of Finite
Elements. Courier Dover Publications, New York (2012)

\bibitem{russo}Russo, A.D., Alonso, A.E.: A posteriori error estimates for nonconforming approximations of Steklov eigenvalue problems.
 Comput. Math. Appl. 62(11), 4100-4117 (2011)

\bibitem{savare}Savar$\acute{e}$, G.: Regularity results for elliptic equations in Lipschitz domains.
J. Funct. Anal. 152, 176-201 (1998)

\bibitem{shiwang} Shi, Z., Wang, M.: Finite Element Methods. Science Press, Beijing (2013)

\bibitem{strang}Strang, G., Fix, G.J.: An alalysis of the finite element method. Prentice-Hall, New York (1973)

\bibitem{verfurth}
Verf$\ddot{u}$rth, R.: A review of a posteriori error
estimates and adaptive mesh-refinement techniques. Wiley-Teubner,
New York (1996)

\bibitem{xie}Xie, H.: A type of multilevel method for the Steklov eigenvalue problem.
IMA J. Numer. Anal. 34, 592-608 (2014)

\bibitem{yang1}Yang, Y., Li, Q., Li, S.: Nonconforming finite element approximations of
the Steklov eigenvalue problem.  Appl. Numer. Math. 59, 2388-2401 (2009)









\end{thebibliography}
\end{document}